# ON PERFECTLY MATCHED LAYERS FOR DISCONTINUOUS PETROV–GALERKIN METHODS

ALI VAZIRI ASTANEH*[1,2], BRENDAN KEITH[1], LESZEK DEMKOWICZ[1]

[1]The Institute for Computational Engineering and Sciences (ICES), The University of Texas at Austin, 201 E 24th St, Austin, TX 78712, USA

[2]MSC Software Corporation, Newport Beach, CA 92660, USA

Abstract. In this article, several discontinuous Petrov–Galerkin (DPG) methods with perfectly matched layers (PMLs) are derived along with their quasi-optimal graph test norms. Ultimately, two different complex coordinate stretching strategies are considered in these derivations. Unlike with classical formulations used by Bubnov–Galerkin methods, with so-called ultraweak variational formulations, these two strategies in fact deliver different formulations in the PML region. One of the strategies, which is argued to be more physically natural, is employed for numerically solving two- and three-dimensional time-harmonic acoustic, elastic, and electromagnetic wave propagation problems, defined in unbounded domains. Through these numerical experiments, efficacy of the new DPG methods with PMLs is verified.

## 1. Introduction

For various reasons, in recent years, there has been a growing interest in non-standard (i.e. non-Bubnov–Galerkin) finite element methods for wave propagation problems. Some of this interest stems from the following features which defy most conventional wisdom from the classical setting. Simultaneously, some such methods: (1) deliver Hermitian *positive-definite* stiffness matrices [1, 2]; (2) offer mesh- and wavenumber-independent discrete stability [3, 4]; (3) yield to desirable iterative solvers, such as conjugate-gradient and multigrid [1, 4–6]; (4) readily permit simple and robust adaptive mesh refinement strategies, even on coarse initial meshes with far under-resolved solutions [4]. Each of the features just described emanate from underlying minimum residual principles which can be formulated in general functional settings [7].

In this article, we will consider wave propagation problems arising from time-harmonic acoustic, electromagnetic, and elastodynamic models. Our focus will be directed towards discontinuous Petrov–Galerkin (DPG) methods which, by their nature, involve such non-symmetric variational formulations [8]. DPG methods have already been studied for acoustic wave equations in [3, 4, 6, 7, 9–14], for Maxwell's equations [15], for elastic media [16–24], and even for applications in nonlinear optics with the Schrödinger equation [25]. Until now, the construction of perfectly matched layers (PMLs) for DPG methods have not been significantly analyzed.

Many wave propagation problems in acoustics, elastodynamics, and electromagnetics are posed on unbounded domains. In computation, this unbounded domain must be partitioned into a bounded interior, wherein the solution is of interest, and an unbounded exterior, wherein the solution can be neglected. By specifying an appropriate absorbing boundary condition (ABC) which mimics the wave absorption properties of the unbounded exterior, all computations can be performed on the bounded

*Corresponding author. *E-mail address*: ali@vaziri.info.

**Acknowledgment.** This work was partially supported with grants by NSF (DMS-1418822), AFOSR (FA9550-12-1-0484), and ONR (N00014-15-1-2496). The first author was also supported in part by the 2016 Peter O'Donnell, Jr. Postdoctoral Fellowship from The Institute for Computational Engineering and Sciences at The University of Texas at Austin. The second author was also supported in part by the 2017 Graduate School University Graduate Continuing Fellowship at The University of Texas at Austin.

interior and the truncated exterior. The most popular ABCs are induced by infinite elements [26, 27], exact non-reflecting boundary conditions [28], rational ABCs and continued-fraction ABCs [29–33], and PMLs [34].

Among these various ABCs, the PML technique has become particularly popular due its accuracy, versatility, and simplicity. The PML was first introduced by Bérenger [34] who showed that dampening out propagating waves without creating artificial reflections (i.e. perfect matching) can be achieved by modifying the conductivity parameters of Maxwell's equations. The complex coordinate stretching approach presented by Chew et al. [35, 36] accelerated subsequent developments due to greater generality and flexibility. Since then, the PML technique has been extensively used for an abundance of wave propagation problems in acoustics [37–43], elastodynamics [42, 44–50], and electromagnetics [35, 42, 51, 52].

One of the challenges in PML applications is stability issues in anisotropic media [53–59] due to admitting wavemodes with differing phase and group velocity signs. Developing well-posed (or stable) PML for such media has been a topic of several studies [46, 53, 56–58]. Devising an efficient PML leading to robust absorption by using the minimal number of layers is also of great importance [60–62]. Many studies have investigated improving the performance of PML by optimizing the PML parameters [62–64], constructing adaptive discretizations [42, 65–67], and using perfectly matched discrete layers (PMDLs) [68–72].

Recognizing these open problems, the focus of this paper is not the stability analysis or optimal discretization of PMLs. Instead, we only lay a simple foundation for their construction in non-symmetric functional settings. Our focus is on so-called (broken) ultraweak variational formulations which are commonly used with DPG methods and have the special feature that they may be used with general polygonal elements [73]. However, our methods can very easily be applied to any other conceivable variational formulation. Complex coordinate stretching has been discussed in an earlier work for different equations using Bubnov–Galerkin methods [74]. Inspired by this and to elaborate the differences, we present two possible categories of PML-consistent ultraweak variational formulations. Then we discuss each category in the framework of DPG methods including the appropriate choice of the adjoint graph norm.

In Section 2, we briefly review the transformation laws in standard exact sequence Hilbert spaces. Then, in Section 3, we summarily introduce DPG methods. In Sections 4–6, we perform complex coordinate stretching and derive ultraweak variational formulations for acoustic, electromagnetic, and elastodynamic wave propagation problems. In Sections 7 and 8, we describe the numerical verification experiments we performed on the respective wave propagation problems. Finally, in Section 9, we summarize our findings.

As a last remark, although we perform numerical verification in both two and three spatial dimensions for the wave propagation problems mentioned above, for overall brevity, we will often only consider the three-dimensional setting in our analysis. We thereby leave the lesser details of the two-dimensional analysis to the reader.

## 2. Preliminaries I: Hilbert spaces & transformation rules

**2.1. The canonical Hilbert spaces.** Let $\Omega \subseteq \mathbb{R}^3$ be a connected domain. Throughout this article, we will often be concerned with energy spaces and operators from the following well-known de Rham complex:

$$H^1(\Omega) \xrightarrow{\nabla} \boldsymbol{H}(\mathrm{curl},\Omega) \xrightarrow{\nabla\times} \boldsymbol{H}(\mathrm{div},\Omega) \xrightarrow{\nabla\cdot} L^2(\Omega)\,. \tag{2.1}$$

In some circumstances, we will require products of two of these energy spaces, $\boldsymbol{H}^1(\Omega) = (\boldsymbol{H}^1(\Omega))^3$ and $\boldsymbol{L}^2(\Omega) = (\boldsymbol{L}^2(\Omega))^3$, or even the symmetric part of the product of the energy space $\boldsymbol{H}(\mathrm{div},\Omega)$:

$$\boldsymbol{H}(\mathrm{div},\Omega;\mathbb{S}) = \left\{\boldsymbol{\sigma} \in (\boldsymbol{H}(\mathrm{div},\Omega))^3 \mid \boldsymbol{\sigma} = \boldsymbol{\sigma}^{\mathsf{T}}\right\}.$$

In these circumstances, $\boldsymbol{H}^1(\Omega)$ is endowed with the row-wise distributional gradient, $\boldsymbol{\nabla}$, and $\boldsymbol{H}(\mathrm{div},\Omega,\mathbb{S})$ with the row-wise distributional divergence, $\boldsymbol{\nabla}\cdot$.

2.2. **Transformation rules in the canonical Hilbert spaces.** Let $\boldsymbol{\varphi}:\mathbb{R}^d\to\mathbb{C}^d$ be a smooth invertible transformation with smooth inverse (a diffeomorphism) $\boldsymbol{x}\mapsto\widetilde{\boldsymbol{x}}=\boldsymbol{\varphi}(\boldsymbol{x})$ and define $\boldsymbol{J}_{ij}=\frac{\partial\widetilde{x}_i}{\partial x_j}$ and $J=\det(\boldsymbol{J})$. The PML variational formulations we will soon derive will come from equations posed on a stretched domain $\widetilde{\Omega}=\boldsymbol{\varphi}(\Omega)$ defined by such a transformation. In order to interpret the corresponding stretched equations on the unstretched computational domain $\Omega$, a coordinate transformation must be applied. In any setting, in order to maintain consistency of the equations, these transformations must be consistent with the energy spaces and operators and, likewise, the entire de Rham complex. Such transformations are called pullbacks in differential geometry and sometimes Piola transformations in engineering literature.

We choose to forgo the details involved in deriving the Piola transformations for the de Rham complex (2.1) because they can be found in various sources. Instead, simply let $\widetilde{u}\in H^1(\widetilde{\Omega})$, $\widetilde{\boldsymbol{E}}\in\boldsymbol{H}(\mathrm{curl},\widetilde{\Omega})$, $\widetilde{\boldsymbol{V}}\in\boldsymbol{H}(\mathrm{div},\widetilde{\Omega})$, and $\widetilde{q}\in L^2(\widetilde{\Omega})$ where

$$H^1(\widetilde{\Omega})\xrightarrow{\widetilde{\nabla}}\boldsymbol{H}(\mathrm{curl},\widetilde{\Omega})\xrightarrow{\widetilde{\nabla}\times}\boldsymbol{H}(\mathrm{div},\widetilde{\Omega})\xrightarrow{\widetilde{\nabla}\cdot}L^2(\widetilde{\Omega})\,.$$

Following [74], the Piola tranformations for the energy spaces above are defined as follows:

$$H^1:\quad \widetilde{u}\circ\boldsymbol{\varphi}=u\,, \tag{2.2a}$$

$$\boldsymbol{H}(\mathrm{curl}):\quad \widetilde{\boldsymbol{E}}\circ\boldsymbol{\varphi}=\boldsymbol{J}^{-\mathsf{T}}\boldsymbol{E}\quad\text{and}\quad\big(\widetilde{\nabla}\widetilde{u}\big)\circ\boldsymbol{\varphi}=\boldsymbol{J}^{-\mathsf{T}}\big(\nabla u\big)\,, \tag{2.2b}$$

$$\boldsymbol{H}(\mathrm{div}):\quad \widetilde{\boldsymbol{V}}\circ\boldsymbol{\varphi}=J^{-1}\boldsymbol{J}\boldsymbol{V}\quad\text{and}\quad\big(\widetilde{\nabla}\times\widetilde{\boldsymbol{E}}\big)\circ\boldsymbol{\varphi}=J^{-1}\boldsymbol{J}\big(\nabla\times\boldsymbol{E}\big)\,, \tag{2.2c}$$

$$L^2:\quad \widetilde{q}\circ\boldsymbol{\varphi}=J^{-1}q\quad\text{and}\quad\big(\widetilde{\nabla}\cdot\widetilde{\boldsymbol{V}}\big)\circ\boldsymbol{\varphi}=J^{-1}\big(\nabla\cdot\boldsymbol{V}\big)\,. \tag{2.2d}$$

Here, $\boldsymbol{J}^{-\mathsf{T}}=(\boldsymbol{J}^{-1})^{\mathsf{T}}$ denotes the inverse and transpose of the Jacobian. Explicitly, $\boldsymbol{J}^{\mathsf{T}}$ is differentiated from the conjugate transpose $\boldsymbol{J}^{\mathsf{H}}=\overline{\boldsymbol{J}}^{\mathsf{T}}$, where $\overline{\boldsymbol{J}}$ is the complex conjugate of $\boldsymbol{J}$. Under the transformations (2.2a)–(2.2d), the corresponding pullbacks belong to the proper energy spaces in (2.1); $u\in H^1(\Omega)$, $\boldsymbol{E}\in\boldsymbol{H}(\mathrm{curl},\Omega)$, $\boldsymbol{V}\in\boldsymbol{H}(\mathrm{div},\Omega)$, and $q\in L^2(\Omega)$.

*Remark* 2.1. Useful Piola transformations also exist for vector fields defined on $\partial\widetilde{\Omega}$. For instance, after accounting for re-normalization,

$$\widetilde{\boldsymbol{n}}\circ\boldsymbol{\varphi}=\frac{J\,\boldsymbol{J}^{-\mathsf{T}}\boldsymbol{n}}{|J\,\boldsymbol{J}^{-\mathsf{T}}\boldsymbol{n}|}\,, \tag{2.3}$$

where $\widetilde{\boldsymbol{n}}$ is normal to $\partial\widetilde{\Omega}$ and $\boldsymbol{n}$ is normal to $\partial\Omega$.

*Remark* 2.2. Using (2.2d), we can compute the pullback of a volume element, $\mathrm{d}\widetilde{\boldsymbol{x}}=J\mathrm{d}\boldsymbol{x}$. Similarly, for boundary integrals, it can be shown that $J_{\partial\Omega}=|J\boldsymbol{J}^{-\mathsf{T}}\boldsymbol{n}|$, where $J_{\partial\Omega}\mathrm{d}S$ is the volume element on $\partial\Omega$.

2.3. **Complex coordinate stretching.** Motivated by [36], we will form PML variational formulations on bounded domains for wave scattering problems on *unbounded domains* by solving for stretched waves (e.g. $\widetilde{u}^\infty\circ\boldsymbol{\varphi}$, where $\widetilde{u}^\infty\in H^1(\widetilde{\Omega})$) instead of their unstretched counterparts (e.g. $u^\infty\in H^1(\Omega)$). The trick then is to construct the corresponding variational formulations in such a way that, in a domain of interest, $\Omega_0$, the stretched solution coincides with its unstretched counterpart (e.g. $\widetilde{u}^\infty\circ\boldsymbol{\varphi}|_{\Omega_0}=u^\infty|_{\Omega_0}$). To sufficiently approximate this embedding without knowing the boundary conditions for the domain of interest, a small layer neighbouring the domain of interest (i.e. the PML) is used and a solution $u^{\mathrm{PML}}|_{\Omega_0}\approx u^\infty|_{\Omega_0}$ is solved for. Here, the stretching function is constructed so that, for outgoing waves, both the stretched outgoing and the stretched incoming waves vanish exponentially. Generally, in order induce such exponential decay, the stretching function maps into in complex plane.

Extensive research has been devoted to designing stretching functions for various problems. In this article, we choose to consider only a relatively simple but important class of stretching functions for Cartesian domains, which we now define.

**Definition 2.1.** Let $L > 0$ be the computational domain length and $0 < l < L$ be the interior domain length. In Cartesian coordinates, define a uniaxial stretch to be any transformation such that:

$$\widetilde{x}_j = \begin{cases} x_j\,, & \text{if } 0 < x_j \leq l\,, \\ x_j + if(x_j,\omega)\,, & \text{if } l < x_j \leq L\,, \end{cases} \tag{2.4}$$

where $i = \sqrt{-1}$ is the imaginary unit. For now, $f(x_j,\omega) > 0$ is arbitrary. Notice that, because $f(x_j,\omega)$ is everywhere positive, the outgoing propagative wavemodes $e^{i\boldsymbol{k}_j\widetilde{\boldsymbol{x}}}$ will decay exponentially in the PML region $l < x_j \leq L$.

*Remark* 2.3. The Jacobian $\boldsymbol{J}$, for a uniaxial stretch in Cartesian coordinates, is diagonal:

$$\boldsymbol{J} = \begin{bmatrix} \frac{\partial\widetilde{x}_1}{\partial x_1} & 0 & 0 \\ 0 & \frac{\partial\widetilde{x}_2}{\partial x_2} & 0 \\ 0 & 0 & \frac{\partial\widetilde{x}_3}{\partial x_3} \end{bmatrix}. \tag{2.5}$$

## 3. Preliminaries II: DPG methods

3.1. **Broken Hilbert spaces and interface spaces.** In the variational formulations to follow, we require several mesh-dependent Hilbert spaces: namely, broken Hilbert spaces and their complementary interface spaces. These are special Hilbert spaces involved in the broken variational formulations which DPG methods generally rely on.

Let $\mathcal{T}$ be a finite Lipschitz partition of $\Omega$, from here on referred to as a mesh. Now define the mesh skeleton $\mathcal{S} = \{\partial K | K \in \mathcal{T}\}$. Roughly following [15], we define the broken Hilbert spaces

$$H^1(\mathcal{T}) = \prod_{K\in\mathcal{T}} H^1(K)\,, \qquad \boldsymbol{H}(\mathrm{curl},\mathcal{T}) = \prod_{K\in\mathcal{T}} \boldsymbol{H}(\mathrm{curl},K)\,, \qquad \text{and} \qquad \boldsymbol{H}(\mathrm{div},\mathcal{T}) = \prod_{K\in\mathcal{T}} \boldsymbol{H}(\mathrm{div},K)\,.$$

Lastly, defining $\boldsymbol{n}_K$ as the normal to an element $K \in \mathcal{T}$, we also require the following interface spaces on the mesh skeleton:

$$H^{1/2}(\mathcal{S}) = \left\{ \widehat{u} \in \prod_{K\in\mathcal{T}} H^{1/2}(\partial K) \;\middle|\; \exists\, u \in H^1(\Omega) \text{ where } u|_{\partial K} = \widehat{u}|_{\partial K}\ \forall K \in \mathcal{T} \right\},$$

$$H^{-1/2}(\mathcal{S}) = \left\{ \widehat{\sigma}_n \in \prod_{K\in\mathcal{T}} H^{-1/2}(\partial K) \;\middle|\; \exists\, \boldsymbol{\sigma} \in \boldsymbol{H}(\mathrm{div},\Omega) \text{ where } (\boldsymbol{\sigma}\cdot\boldsymbol{n}_K)|_{\partial K} = \widehat{\sigma}_n|_{\partial K}\ \forall K \in \mathcal{T} \right\},$$

and

$$\boldsymbol{H}^{-1/2}(\mathrm{curl},\mathcal{S}) = \Bigg\{ \widehat{\mathbf{E}} \in \prod_{K\in\mathcal{T}} \boldsymbol{H}^{-1/2}(\mathrm{curl},\partial K) \;\Big|\; \exists\, \mathbf{E} \in \boldsymbol{H}(\mathrm{curl},\Omega) \text{ where } ((\boldsymbol{n}_K \times \mathbf{E}) \times \boldsymbol{n}_K)|_{\partial K} = \widehat{\mathbf{E}}|_{\partial K}\ \forall K \in \mathcal{T} \Bigg\},$$

Similar definitions as above hold for products of these spaces. e.g. $\boldsymbol{H}^1(\mathcal{T}) = \prod_{K\in\mathcal{T}} \boldsymbol{H}^1(K)$.

Many important properties of the above spaces are known but are not necessary for our analysis. We refer the interested reader to the detailed discussion in [15, Section 2].

3.2. **Broken variational formulations.** Consider the following frequency-domain continuity equation for a given quantity, say $q$:

$$-\,i\omega\rho - \nabla\cdot\mathbf{j} = g\,.$$

Here, $\omega$ is the temporal frequency, $\rho$ is the volume density of $q$, $\mathbf{j}$ is the flux of $q$, and $g$ is a source generating $q$. In the next three sections, we will encounter several equations like this. Formally, multiplying the continuity equation equation by a test function, $v$, integrating over a single element, $K \in \mathcal{T}$, and then integrating by parts, we obtain

$$-\,i\omega(\rho,v)_K + (\mathbf{j},\nabla v)_K - \langle \mathbf{j}\cdot\boldsymbol{n}_K, v\rangle_{\partial K} = (g,v)_K\,.$$

Here, $(\cdot,\cdot)_K$ is the standard $L^2(K)$-inner product and $\langle\cdot,\cdot\rangle_{\partial K}$ is a duality pairing; assuming sufficiently smooth arguments, it is also the $L^2(\partial K)$-inner product.

In a common formal derivation of broken variational formulations, the flux variable $\mathbf{j}$ is distinguished on the element $K$ and on its boundary $\partial K$ by introducing a new interface flux variable, $\widehat{\mathrm{j}}_n$, to replace $\mathbf{j}\cdot\boldsymbol{n}_K$. Now, summing over each element in the mesh, we arrive at the equation

$$-i\omega(\rho,v)_\Omega + (\mathbf{j},\nabla_h v)_\Omega - \langle\widehat{\mathrm{j}}_n, v\rangle_{\partial\mathcal{T}} = (g,v)_\Omega\,. \tag{3.1}$$

Here, $\nabla_h$ denotes the *element-wise* distributional gradient and $\langle\widehat{\mathrm{j}}_n, v\rangle_{\partial\mathcal{T}} = \sum_{K\in\mathcal{T}}\langle\widehat{\mathrm{j}}_n, v\rangle_{\partial K}$ denotes the accumulation of all respective boundary terms.

Many variational formulations derived like (3.1), can be proven to be well-posed when the test functions are taken from a broken Hilbert space, $v\in\mathcal{V}(\mathcal{T})$ [15, 21, 75]. In such scenarios, the interface functions naturally belong to interface spaces, $\widehat{\mathrm{j}}_n\in\widehat{\mathcal{U}}(\mathcal{S})$. Interface functions in broken variational formulations act like Lagrange multipliers in response to the discontinuity of the test space.

In this paper, we will consider (broken) ultraweak variational formulations. In these formulations, all of the derivatives have moved onto the test function (e.g. onto $v$, above) through integration by parts.

3.3. **DPG methods.** Among many other sources, implementation of (practical) DPG methods is described in significant detail in [2, 7, 76, 77]. Therefore, we will only briefly describe their key properties.

Let $\mathcal{U}$ and $\mathcal{V}$ be Hilbert spaces, *where* $\mathcal{V}$ *is broken*, and let $\mathcal{U}_h\subseteq\mathcal{U}$ and $\mathcal{V}_r\subseteq\mathcal{V}$ be finite dimensional. Let $\ell\in\mathcal{V}'$ and let $b:\mathcal{U}\times\mathcal{V}\to\mathbb{C}$ be a continuous sesquilinear form. Assuming the appropriate inf–sup condition for $b$ and the existence of Fortin operator $\Pi:\mathcal{V}\to\mathcal{V}_r$ (see [19, 78]), a practical DPG method with optimal test functions is a discrete least-squares finite element method [7] defined by the problem

$$\begin{cases}\text{Find } \mathfrak{u}_h\in\mathcal{U}_h :\\ b(\mathfrak{u}_h,\Theta_r\mathfrak{w}_h) = \ell(\Theta_r\mathfrak{w}_h)\,, \qquad \forall\,\mathfrak{w}_h\in\mathcal{U}_h\,,\end{cases}$$

where $\Theta_r:\mathcal{U}_h\to\mathcal{V}_r$ is the discrete trial-to-test operator defined through the inner product $(\cdot,\cdot)_\mathcal{V}$ from the test space $\mathcal{V}$:

$$(\Theta_r\mathfrak{w}_h,\mathfrak{v}_r)_\mathcal{V} = b(\mathfrak{w}_h,\mathfrak{v}_r)\,,\quad \forall\,\mathfrak{w}_h\in\mathcal{U}_h\;\mathfrak{v}_r\in\mathcal{V}_r\,. \tag{3.2}$$

For general test norms, computation of the trial-to-test operator $\Theta_r$ in (3.2) requires the inversion of a large symmetric Gram matrix corresponding to the discretization of the inner product $(\cdot,\cdot)_\mathcal{V}$. This is made feasible by using a broken test space and special so-called "localizable" test norms, $\|\cdot\|_\mathcal{V}$ [2, 7]; i.e. norms on the entire space $\mathcal{V}$, which, for any Lipschitz $K\in\mathcal{T}$, induce norms, $\|\cdot\|_{\mathcal{V}_K} = \|\cdot|_K\|_\mathcal{V}$, when restricted to $\mathcal{V}_K = \{\mathfrak{v}|_K \mid \mathfrak{v}\in\mathcal{V}\}$. In this case, the Gram matrix becomes block-diagonal and the inverse can be computed locally and in parallel. This is the reason broken variational formulations are essential in DPG methods. The following final preliminary subsection introduces a particularly well-studied localizable test norm often used in DPG methods.

3.4. **The graph norm.** With broken ultraweak variational formulations, it is often the case that the bilinear form can be defined

$$b(\mathfrak{u},\mathfrak{v}) = (\mathrm{u},\mathcal{L}^*\mathfrak{v})_\Omega + \langle\widehat{\mathrm{u}},\mathfrak{v}\rangle_{\partial\mathcal{T}}\,,\quad \forall\,\mathfrak{u}=(\mathrm{u},\widehat{\mathrm{u}})\in\mathcal{U}=L^2(\Omega)\times\widehat{\mathcal{U}},\;\mathfrak{v}\in\mathcal{V}=\mathrm{Dom}(\mathcal{L}^*)\,, \tag{3.3}$$

where $\mathcal{L}^*:\mathcal{V}\to L^2(\Omega)$ be a closed linear operator and $\mathrm{Dom}(\mathcal{L}^*)\subseteq L^2(\Omega)$ is dense. Notice that (3.1) happens to suggest this functional setting.

For various reasons (see [79]), for any continuously inf–sup stable bilinear form $b$, it is desirable to use test norms $\|\cdot\|_\mathcal{V}$ closely related to the so-called *optimal test norm*:

$$|||\mathfrak{v}|||_\mathcal{V} = \sup_{\mathfrak{u}\in\mathcal{U}}\frac{|b(\mathfrak{u},\mathfrak{v})|}{\|\mathfrak{u}\|_\mathcal{U}}\,, \tag{3.4}$$

where $\|\cdot\|_\mathcal{U}$ is a desired norm on the trial space $\mathcal{U}$.

If the test functions in (3.3) are not broken, then the second term in the bilinear form $b$ can be discarded and, simply, $b(\mathfrak{u},\mathfrak{v}) = (\mathfrak{u},\mathcal{L}^*\mathfrak{v})_\Omega$. In this unbroken setting, which is usually impractical for computation, observe that the optimal test norm (3.4) is readily computable for all $\mathfrak{v}\in\mathcal{V}$:

$$|||\mathfrak{v}|||_{\mathcal{V}} = \sup_{\mathfrak{u}\in L^2(\Omega)} \frac{(\mathfrak{u},\mathcal{L}^*\mathfrak{v})_\Omega}{\|\mathfrak{u}\|_\Omega} = \|\mathcal{L}^*\mathfrak{v}\|_\Omega\,.$$

Unfortunately, although this norm is closely related to the optimal test norm for the broken test space setting, this norm is not localizable in general. Fortunately, it is usually equivalent to the so-called (*adjoint*) *graph norm*:

$$\|\mathfrak{v}\|^2_{\mathcal{L}^*} = \|\mathcal{L}^*\mathfrak{v}\|^2_\Omega + \|\mathfrak{v}\|^2_\Omega\,,$$

which is localizable [3]. This is the test norm we will consider in this article.

## 4. Acoustics

When deriving a variational formulation for a model with PMLs, there are two natural paths to consider:

(1) Beginning with the strong form of the equations in the stretched coordinates, *first* construct the pulled-back equations, *then* multiply with test functions and integrate by parts in the spatial coordinates.
(2) Alternatively, *first* multiply with test functions and integrate by parts to construct a variational formulation in the stretched coordinates, *then* pull back the corresponding stretched variational formulation to spatial coordinates.

In classical functional settings, it is well-known that these two approaches are identical in that they return the same ultimate variational formulation [74]. However, there are subtle differences to acknowledge in alternative functional settings. We will now illustrate these differences by example; namely, deriving the acoustic wave equations in the ultraweak setting. Other common settings [15] can be treated similarly.

### 4.1. A variational formulation from the pulled-back equations.

Here, we contemplate the first approach described above. Using the time-harmonic ansatz $\widetilde{p}(\widetilde{\boldsymbol{x}},t) = \widetilde{p}(\widetilde{\boldsymbol{x}})e^{-i\omega t}$, consider the strong form of the acoustic wave equations with angular frequency $\omega$:

$$\left\{\begin{aligned} -i\omega\widetilde{\mathbf{u}} - \widetilde{\nabla}\widetilde{p} &= 0\,,\\ -i\omega\widetilde{p} - \widetilde{\nabla}\cdot\widetilde{\mathbf{u}} &= i\omega^{-1}\widetilde{f}\,, \end{aligned}\right. \iff \left\{\begin{aligned} -i\omega J^{-1}\boldsymbol{J}\mathbf{u} - \boldsymbol{J}^{-\mathsf{T}}\nabla p &= 0\,,\\ -i\omega p - J^{-1}\nabla\cdot\mathbf{u} &= i\omega^{-1}J^{-1}f\,. \end{aligned}\right. \tag{4.1}$$

On the left-hand side are the equations in the stretched coordinates; $\widetilde{\mathbf{u}}$ is the velocity and $\widetilde{p}$ is the pressure. On the right-hand side are the corresponding pulled-back equations in the spatial coordinates; $\mathbf{u}$ is the pulled-back velocity and $p$ is the pulled-back pressure. Note that we have chosen to assume $\widetilde{f}\in L^2(\widetilde{\Omega})$ and that it obeys the transformation rule (2.2d).

After pre-multiplying the first and second pulled-back equation by $\boldsymbol{J}^{\mathsf{T}}$ and $J$, respectively, it is easy to construct the corresponding variational formulation. Indeed, after performing this pre-multiplication, multiply the pulled-back equations by test functions $\mathbf{v}$ and $q$, respectively, and integrate each equation over the domain $\Omega$:

$$\left\{\begin{aligned} (-i\omega J^{-1}\boldsymbol{J}^{\mathsf{T}}\boldsymbol{J}\mathbf{u},\mathbf{v})_\Omega - (\nabla p,\mathbf{v})_\Omega &= 0\,,\\ (-i\omega Jp,q)_\Omega - (\nabla\cdot\mathbf{u},q)_\Omega &= (i\omega^{-1}f,q)_\Omega\,. \end{aligned}\right. \tag{4.2}$$

Notice that element-wise integration by parts can now be performed freely in both equations above, after which, we arrive at the following variational formulation:

$$\left\{\begin{aligned} &\text{Find } (p,\mathbf{u},\widehat{p},\widehat{\mathrm{u}}_n)\in L^2(\Omega)\times\boldsymbol{L}^2(\Omega)\times H^{1/2}_0(\mathcal{S})\times H^{-1/2}(\mathcal{S}):\\ &(-i\omega J^{-1}\boldsymbol{J}^{\mathsf{T}}\boldsymbol{J}\mathbf{u},\mathbf{v})_\Omega + (p,\nabla_h\cdot\mathbf{v})_\Omega - \langle\widehat{p},\mathbf{v}\cdot\boldsymbol{n}\rangle_{\partial\mathcal{T}} = 0\,, &&\forall\,\mathbf{v}\in\boldsymbol{H}(\mathrm{div},\mathcal{T})\,,\\ &\qquad(-i\omega Jp,q)_\Omega + (\mathbf{u},\nabla_h q)_\Omega - \langle\widehat{\mathrm{u}}_n,q\rangle_{\partial\mathcal{T}} = (i\omega^{-1}f,q)_\Omega\,, &&\forall\,q\in H^1(\mathcal{T})\,. \end{aligned}\right. \tag{4.3}$$

Here, we have chosen to include the corresponding functional settings, which are valid if $\boldsymbol{J}, \boldsymbol{J}^{-1} \in L^\infty(\Omega)$, and also homogeneous boundary conditions on the pressure.

*Remark* 4.1. Denote $\mathfrak{v} = (\mathbf{v}, q)$ and define $\mathcal{V} = \boldsymbol{H}(\mathrm{div}, \mathcal{T}) \times H^1(\mathcal{T})$. The corresponding adjoint graph norm is simply

$$\|\mathfrak{v}\|_{\mathcal{L}^*}^2 = \|i\omega(J^{-1}\boldsymbol{J}^\mathsf{T}\boldsymbol{J})^\mathsf{H}\mathbf{v} + \nabla_h q\|_\Omega^2 + \|i\omega\overline{J}q + \nabla_h \cdot \mathbf{v}\|_\Omega^2 + \|q\|_\Omega^2 + \|\mathbf{v}\|_\Omega^2\,, \quad \forall\, \mathfrak{v} \in \mathcal{V}\,.$$

4.2. **Pulling back a stretched variational formulation.** Here, we contemplate the alternative approach described above. Beginning with the strong formulation of acoustics in stretched coordinates (4.1), multiply with arbitrary test functions $\widetilde{\mathbf{v}}$ and $\widetilde{q}$, and immediately integrating by parts in the stretched coordinates. Now, observe that we arrive at the following broken variational formulation:

$$(4.4)\qquad \begin{cases} \text{Find } (\widetilde{p}, \widetilde{\mathbf{u}}, \widetilde{\widehat{p}}, \widetilde{\widehat{\mathrm{u}}}_{\widetilde{n}}) \in L^2(\widetilde{\Omega}) \times \boldsymbol{L}^2(\widetilde{\Omega}) \times H_0^{1/2}(\widetilde{\mathcal{S}}) \times H^{-1/2}(\widetilde{\mathcal{S}}): \\ (-i\omega\widetilde{\mathbf{u}}, \widetilde{\mathbf{v}})_{\widetilde{\Omega}} + (\widetilde{p}, \widetilde{\nabla_h} \cdot \widetilde{\mathbf{v}})_{\widetilde{\Omega}} - \langle \widetilde{\widehat{p}}, \widetilde{\mathbf{v}} \cdot \widetilde{\boldsymbol{n}} \rangle_{\partial\widetilde{\mathcal{T}}} = 0\,, & \forall\, \widetilde{\mathbf{v}} \in \boldsymbol{H}(\mathrm{div}, \widetilde{\mathcal{T}})\,, \\ \qquad (-i\omega\widetilde{p}, \widetilde{q})_{\widetilde{\Omega}} + (\widetilde{\mathbf{u}}, \widetilde{\nabla_h}\widetilde{q})_{\widetilde{\Omega}} + \langle \widetilde{\widehat{\mathrm{u}}}_{\widetilde{n}}, \widetilde{q} \rangle_{\partial\widetilde{\mathcal{T}}} = (i\omega^{-1}\widetilde{f}, \widetilde{q})_{\widetilde{\Omega}}\,, & \forall\, \widetilde{q} \in H^1(\widetilde{\mathcal{T}})\,. \end{cases}$$

Again, homogeneous boundary conditions have been imposed on the pressure.

Let us fix the test function $(\widetilde{\mathbf{v}}, \widetilde{q}) \in \boldsymbol{H}(\mathrm{div}, \widetilde{\mathcal{T}}) \times H^1(\widetilde{\mathcal{T}})$ in order to construct the pullback of the stretched variational formulation (4.4). For the sake of rigor, we choose to let $\widetilde{p}_{\mathrm{ext}}$ be an $H^1$-extension of $\widetilde{\widehat{p}}$ and, likewise, let $\widetilde{\mathbf{u}}_{\mathrm{ext}}$ be an $\boldsymbol{H}(\mathrm{div})$-extension of $\widetilde{\widehat{\mathrm{u}}}_{\widetilde{n}}$. Therefore, $\mathrm{tr}_{\mathrm{grad}}^{\partial\widetilde{\mathcal{T}}}\,\widetilde{p}_{\mathrm{ext}} = \widetilde{\widehat{p}}$ and $\mathrm{tr}_{\mathrm{div}}^{\partial\widetilde{\mathcal{T}}}\,\widetilde{\mathbf{u}}_{\mathrm{ext}} = \widetilde{\widehat{\mathrm{u}}}_{\widetilde{n}}$.

Now, restricted to a single element $K \in \mathcal{T}$, the first variational equation becomes

$$(-i\omega J(J^{-1}\mathbf{u}), J^{-1}\boldsymbol{J}\mathbf{v})_K + (J(J^{-1}p), J^{-1}\nabla \cdot \mathbf{v})_K - \Big\langle J_{\partial K} p_{\mathrm{ext}}, J^{-1}\boldsymbol{J}\mathbf{v} \cdot \frac{\boldsymbol{J}\boldsymbol{J}^{-\mathsf{T}}\boldsymbol{n}}{|\boldsymbol{J}\boldsymbol{J}^{-\mathsf{T}}\boldsymbol{n}|} \Big\rangle_{\partial K} = 0\,.$$

where $p_{\mathrm{ext}}$ is the $H^1$-pullback of $\widetilde{p}_{\mathrm{ext}}$. Meanwhile, the second variational equation becomes

$$(-i\omega J(J^{-1}p), q)_K + (J(J^{-1}\mathbf{u}), \boldsymbol{J}^{-\mathsf{T}}\nabla q)_K - \Big\langle J_{\partial K}(J^{-1}\boldsymbol{J}\mathbf{u}_{\mathrm{ext}}) \cdot \frac{\boldsymbol{J}\boldsymbol{J}^{-\mathsf{T}}\boldsymbol{n}}{|\boldsymbol{J}\boldsymbol{J}^{-\mathsf{T}}\boldsymbol{n}|}, q \Big\rangle_{\partial K} = (i\omega^{-1}f, q)_K\,,$$

where $\mathbf{u}_{\mathrm{ext}}$ is the $\boldsymbol{H}(\mathrm{div})$-pullback of $\widetilde{\mathbf{u}}_{\mathrm{ext}}$. Sum over over all elements and recall that $J_{\partial K} = |\boldsymbol{J}\boldsymbol{J}^{-\mathsf{T}}\boldsymbol{n}|$. Denoting $\widehat{p} = \mathrm{tr}_{\mathrm{grad}}^{\partial\mathcal{T}}\,p_{\mathrm{ext}}$ and $\widehat{\mathrm{u}}_n = \mathrm{tr}_{\mathrm{div}}^{\partial\mathcal{T}}\,\mathbf{u}_{\mathrm{ext}}$, the two pulled-back variational equations above immediately simplify to

$$(4.5)\qquad \begin{cases} \text{Find } (p, \mathbf{u}, \widehat{p}, \widehat{\mathrm{u}}_n) \in L^2(\Omega) \times \boldsymbol{L}^2(\Omega) \times H_0^{1/2}(\mathcal{S}) \times H^{-1/2}(\mathcal{S}): \\ (-i\omega\mathbf{u}, J^{-1}\boldsymbol{J}\mathbf{v})_\Omega + (p, J^{-1}\nabla_h \cdot \mathbf{v})_\Omega - \langle \widehat{p}, \mathbf{v} \cdot \boldsymbol{n} \rangle_{\partial\mathcal{T}} = 0\,, & \forall\, \mathbf{v} \in \boldsymbol{H}(\mathrm{div}, \mathcal{T})\,, \\ \qquad (-i\omega p, q)_\Omega + (\mathbf{u}, \boldsymbol{J}^{-\mathsf{T}}\nabla_h q)_\Omega - \langle \widehat{\mathrm{u}}_n, q \rangle_{\partial\mathcal{T}} = (i\omega^{-1}f, q)_\Omega\,, & \forall\, q \in H^1(\mathcal{T})\,. \end{cases}$$

Here, we have again introduced the natural functional setting, assuming $\boldsymbol{J}, \boldsymbol{J}^{-1} \in L^\infty(\Omega)$.

*Remark* 4.2. As before, denote $\mathfrak{v} = (\mathbf{v}, q)$ and define $\mathcal{V} = \boldsymbol{H}(\mathrm{div}, \mathcal{T}) \times H^1(\mathcal{T})$. The corresponding graph norm is now simply

$$(4.6)\qquad \|\mathfrak{v}\|_{\mathcal{L}^*}^2 = \|i\omega J^{-1}\boldsymbol{J}\mathbf{v} + \boldsymbol{J}^{-\mathsf{T}}\nabla_h q\|_\Omega^2 + \|i\omega q + J^{-1}\nabla_h \cdot \mathbf{v}\|_\Omega^2 + \|q\|_\Omega^2 + \|\mathbf{v}\|_\Omega^2\,, \quad \forall\, \mathfrak{v} \in \mathcal{V}\,.$$

*Remark* 4.3. It is natural to consider whether the adjoint graph norm corresponding to the stretched variational formulation (4.4) will be mapped to the adjoint graph norm (4.6). A simple calculation shows that this is not the case.

4.3. **Comparison of the approaches.** Clearly, variational formulations (4.3) and (4.5) are not the same. Nevertheless, both are viable for computation and we will not provide a definitive argument to use one over the other.

Note that the converged values of the interface unknowns $\widehat{\mathrm{u}}_n$ and $\widehat{p}$ will different between the two equations, since they are themselves Lagrange multipliers determined by the corresponding *test space*. In the derivation of the first variational formulation (4.3), it should already be clear at the point of (4.2) that the test functions naturally come from the physical domain and it is the trial functions which have been pulled-back to the physical domain. Therefore, upon integration by parts, the interface unknowns can be identified with traces of the pulled-back solution $\mathrm{tr}^{\partial\mathcal{T}}_{\mathrm{grad}}\, p \sim \widehat{p}$ and $\mathrm{tr}^{\partial\mathcal{T}}_{\mathrm{div}}\, \mathbf{u} \sim \widehat{\mathrm{u}}_n$ (see [15] for pertinent discussion).

Alternatively, consider (4.5). Using the transformation rules (2.2a)–(2.2d) and $\mathcal{V} \xrightarrow{\varphi} \widetilde{\mathcal{V}}$, note that it may be more naturally expressed as

$$\begin{cases} \text{Find } (p, \mathbf{u}, \widehat{p}, \widehat{\mathrm{u}}_{\widetilde{n}}) \in L^2(\Omega) \times \boldsymbol{L}^2(\Omega) \times H_0^{1/2}(\mathcal{S}) \times H^{-1/2}(\mathcal{S}) : \\ (-i\omega\mathbf{u}, \widetilde{\mathbf{v}} \circ \boldsymbol{\varphi})_\Omega + (p, (\widetilde{\nabla}_h \cdot \widetilde{\mathbf{v}}) \circ \boldsymbol{\varphi})_\Omega - \langle J_{\partial\mathcal{T}}\widehat{p}, (\widetilde{\mathbf{v}} \cdot \widetilde{\boldsymbol{n}}) \circ \boldsymbol{\varphi}\rangle_{\partial\mathcal{T}} = 0\,, & \forall \widetilde{\mathbf{v}} \in \boldsymbol{H}(\mathrm{div}, \widetilde{\mathcal{T}})\,, \\ \big(-i\omega p, \widetilde{q} \circ \boldsymbol{\varphi}\big)_\Omega + \big(\mathbf{u}, (\widetilde{\nabla}_h \widetilde{q}) \circ \boldsymbol{\varphi}\big)_\Omega - \langle J_{\partial\mathcal{T}}\widehat{\mathrm{u}}_{\widetilde{n}}, \widetilde{q} \circ \boldsymbol{\varphi}\rangle_{\partial\mathcal{T}} = \big(i\omega^{-1} f, \widetilde{q} \circ \boldsymbol{\varphi}\big)_\Omega\,, & \forall \widetilde{q} \in H^1(\widetilde{\mathcal{T}})\,. \end{cases}$$

In this form, the interpretation as an acoustic wave equation with stretched test functions is more clear. Indeed, first consider the conforming (unbroken) test function scenario. Here, the interface terms disappear and we are to only left find $p$ and $\mathbf{u}$, which are clearly seen to be integrated against stretched test functions. To then make sense of the interface terms, we first note that the identification with traces of the pulled-back solution is gone. Indeed, the interface variables in (4.5), $\widehat{p}$ and $\widehat{\mathrm{u}}_n$, act as Lagrangian multipliers on the stretched mesh instead.

Because we (the authors) find comfort in having a physical interpretation of the interface variables, $\widehat{p}$ and $\widehat{\mathrm{u}}_n$, we have chosen to focus only on the first class of PML formulation in the following examples. This is not to say that the other formulations cannot be used successfully. Indeed such variational formulations were analyzed for acoustics and elastodynamics in [18].

*Remark* 4.4. In the dual setting, such as considered by DPG* methods [80, 81], the alternative variational formulation (4.5) may be a more "natural" fit. Mainly, this is because the solution is then found in the space $\mathcal{V}$, instead of $\mathcal{U}$, and the physical interpretation of the interface variable is no longer relevant.

## 5. Electromagnetics

Using the time-harmonic ansatz $\widetilde{\mathbf{E}}(\widetilde{\boldsymbol{x}}, t) = \widetilde{\mathbf{E}}(\widetilde{\boldsymbol{x}})e^{-i\omega t}$, we consider the strong form of the three-dimensional, linear, isotropic Maxwell's equations:

$$\text{(5.1)}\quad \begin{cases} \widetilde{\nabla} \times \widetilde{\mathbf{E}} - i\omega\mu\widetilde{\mathbf{H}} = 0\,, \\ \widetilde{\nabla} \times \widetilde{\mathbf{H}} + (i\omega\epsilon - \sigma)\widetilde{\mathbf{E}} = \widetilde{\mathbf{J}}^{\mathrm{imp}}\,, \end{cases} \iff \begin{cases} J^{-1}\boldsymbol{J}\nabla \times \mathbf{E} - i\omega\mu\boldsymbol{J}^{-\mathsf{T}}\mathbf{H} = 0\,, \\ J^{-1}\boldsymbol{J}\nabla \times \mathbf{H} + (i\omega\epsilon - \sigma)\boldsymbol{J}^{-\mathsf{T}}\mathbf{E} = J^{-1}\boldsymbol{J}\mathbf{J}^{\mathrm{imp}}\,. \end{cases}$$

Here, $\omega$ is the angular frequency, $\epsilon$ is the electric permittivity, $\mu$ is the magnetic permeability, and $\sigma$ is the material conductivity. On the left-hand side are the equations in the stretched coordinates; $\widetilde{\mathbf{E}}$ is the electric field and $\widetilde{\mathbf{H}}$ is the magnetic field. On the right-hand side are the corresponding pulled-back equations in the spatial coordinates; $\mathbf{E}$ is the pulled-back electric field and $\mathbf{H}$ is the pulled-back magnetic field. Note that we have chosen to assume that the impressed current $\widetilde{\mathbf{J}}^{\mathrm{imp}} \in \boldsymbol{H}(\mathrm{div}, \widetilde{\Omega})$. This is physically motivated by the charge conservation equation in a non-conducting material:

$$\widetilde{\nabla} \cdot \widetilde{\mathbf{J}}^{\mathrm{imp}} = i\omega\widetilde{\rho}\,,$$

where $\widetilde{\rho}$ is the volume charge density in the stretched coordinates. It is also mathematically natural when considering the range of the operators at hand (2.1). Ultimately, this naturally implies that $\widetilde{\mathbf{J}}^{\mathrm{imp}}$ undergoes the transformation (2.2c).

Similar to the treatment of the acoustic wave equations in Section 4.1, pre-multiply both the first and second pulled-back equation by $J\boldsymbol{J}^{-1}$. Now upon multiplication with test functions $\mathbf{F}$ and $\mathbf{G}$, respectively, and integration over the computational domain $\Omega$, we arrive at

$$\begin{cases} \big(\nabla \times \mathbf{E}, \mathbf{F}\big)_\Omega - \big(i\omega\mu J\boldsymbol{J}^{-1}\boldsymbol{J}^{-\mathsf{T}}\mathbf{H}, \mathbf{F}\big)_\Omega = 0\,, \\ \big(\nabla \times \mathbf{H}, \mathbf{G}\big)_\Omega + \big((i\omega\epsilon - \sigma)J\boldsymbol{J}^{-1}\boldsymbol{J}^{-\mathsf{T}}\mathbf{E}, \mathbf{G}\big)_\Omega = \big(\mathbf{J}^{\mathrm{imp}}, \mathbf{G}\big)_\Omega\,. \end{cases}$$

From here, element-wise integration by parts can be freely performed and we uncover the following broken ultraweak variational formulation:

$$
(5.2)\quad \begin{cases} \text{Find } (\mathbf{E},\mathbf{H},\widehat{\mathbf{E}},\widehat{\mathbf{H}}) \in \boldsymbol{L}^2(\Omega)\times\boldsymbol{L}^2(\Omega)\times\boldsymbol{H}_0^{-1/2}(\mathrm{curl},\mathcal{S})\times\boldsymbol{H}^{-1/2}(\mathrm{curl},\mathcal{S}) : \\ \qquad \big(\mathbf{E},\nabla_h\times\mathbf{F}\big)_\Omega - \big(i\omega\mu J\boldsymbol{J}^{-1}\boldsymbol{J}^{-\mathsf{T}}\mathbf{H},\mathbf{F}\big)_\Omega + \big\langle\widehat{\mathbf{E}},\boldsymbol{n}\times\mathbf{F}\big\rangle_{\partial\mathcal{T}} = 0\,, \qquad \forall\,\mathbf{F}\in\boldsymbol{H}(\mathrm{curl},\mathcal{T})\,, \\ \big(\mathbf{H},\nabla_h\times\mathbf{G}\big)_\Omega + \big((i\omega\epsilon-\sigma)J\boldsymbol{J}^{-1}\boldsymbol{J}^{-\mathsf{T}}\mathbf{E},\mathbf{G}\big)_\Omega + \big\langle\widehat{\mathbf{H}},\boldsymbol{n}\times\mathbf{G}\big\rangle_{\partial\mathcal{T}} = \big(\mathbf{J}^{\mathrm{imp}},\mathbf{G}\big)_\Omega\,, \\ \qquad\qquad\qquad\qquad\qquad\qquad\qquad\qquad \forall\,\mathbf{G}\in\boldsymbol{H}(\mathrm{curl},\mathcal{T})\,. \end{cases}
$$

Here, we have tacitly closed the problem with perfect electric conductor (PEC) boundary conditions and assigned of the correct functional setting for well-posedness.

*Remark* 5.1. Denote $\mathfrak{v} = (\mathbf{F},\mathbf{G})$ and define $\mathcal{V} = \boldsymbol{H}(\mathrm{curl},\mathcal{T})\times\boldsymbol{H}(\mathrm{curl},\mathcal{T})$. The corresponding adjoint graph norm is simply

$$
\|\mathfrak{v}\|^2_{\mathcal{L}^*} = \|\nabla_h\times\mathbf{F}-(i\omega\epsilon+\sigma)(J\boldsymbol{J}^{-1}\boldsymbol{J}^{-\mathsf{T}})^{\mathsf{H}}\mathbf{G}\|^2_\Omega+\|\nabla_h\times\mathbf{G}+i\omega\mu(J\boldsymbol{J}^{-1}\boldsymbol{J}^{-\mathsf{T}})^{\mathsf{H}}\mathbf{F}\|^2_\Omega+\|\mathbf{F}\|^2_\Omega+\|\mathbf{G}\|^2_\Omega\,, \quad \forall\,\mathfrak{v}\in\mathcal{V}\,.
$$

## 6. Elastodynamics

For simplicity, in this section we only consider the symmetric case $\boldsymbol{J} = \boldsymbol{J}^{\mathsf{T}}$. Notably, the following derivation indeed handles uniaxial stretching (2.5).

Consider the strong form of the homogeneous, time-harmonic linearized elasticity equations:

$$
(6.1)\quad \begin{cases} \mathsf{S}:\widetilde{\boldsymbol{\sigma}}-\widetilde{\boldsymbol{\varepsilon}}(\widetilde{\boldsymbol{u}}) = 0\,, \\ -\widetilde{\boldsymbol{\nabla}}\cdot\widetilde{\boldsymbol{\sigma}}-\rho\omega^2\widetilde{\boldsymbol{u}} = \widetilde{\boldsymbol{f}}\,, \end{cases} \iff \begin{cases} J^{-1}\mathsf{S}:\boldsymbol{J}\boldsymbol{\sigma}-\boldsymbol{J}^{-1}\boldsymbol{\varepsilon}(\boldsymbol{u}) = 0\,, \\ -J^{-1}\,\boldsymbol{\nabla}\cdot\boldsymbol{\sigma}-\rho\omega^2\boldsymbol{u} = J^{-1}\boldsymbol{f}\,. \end{cases}
$$

Here, $\omega$ is the angular frequency, $\mathsf{S}$ is the compliance tensor, and $\rho$ is the (uniform) mass density. Meanwhile, $\widetilde{\boldsymbol{\varepsilon}}(\widetilde{\boldsymbol{u}})$ and $\boldsymbol{\varepsilon}(\boldsymbol{u})$ are the linearized strains, defined through the symmetric gradient operators $\widetilde{\boldsymbol{\varepsilon}}(\cdot) = \frac{1}{2}(\widetilde{\boldsymbol{\nabla}}\,\cdot+\widetilde{\boldsymbol{\nabla}}^{\mathsf{T}}\cdot)$ and $\boldsymbol{\varepsilon}(\cdot) = \frac{1}{2}(\boldsymbol{\nabla}\,\cdot+\boldsymbol{\nabla}^{\mathsf{T}}\cdot)$, respectively. Notice that we have place the body force $\widetilde{\boldsymbol{f}}\in\boldsymbol{L}^2(\widetilde{\Omega})$ and so, component by component, it naturally obeys the transformation rule (2.2d).

As with the acoustic wave equations in Section 4.1, pre-multiply both the first and second pulled-back equation by $\boldsymbol{J}$ and $J$, respectively. Then, multiplying each equation with *symmetric* $\boldsymbol{\tau}$ — that is, $\boldsymbol{\tau}(\boldsymbol{x})\in\mathbb{S}$, at each spatial coordinate $\boldsymbol{x}\in\Omega$ — and the vector field $\boldsymbol{v}$, respectively, and then integrating parts, we find

$$
\begin{cases} \big(J^{-1}\mathsf{S}:\boldsymbol{J}\boldsymbol{\sigma},\overline{\boldsymbol{J}}\boldsymbol{\tau}\big)_\Omega - \big(\boldsymbol{\nabla}\,\boldsymbol{u},\boldsymbol{\tau}\big)_\Omega = 0\,, \\ \quad -\big(\boldsymbol{\nabla}\cdot\boldsymbol{\sigma},\boldsymbol{v}\big)_\Omega - \big(J\rho\omega^2\boldsymbol{u},\boldsymbol{v}\big)_\Omega = \big(\boldsymbol{f},\boldsymbol{v}\big)_\Omega\,. \end{cases}
$$

At this point, after assignment of the correct functional setting for well-posedness, simple element-wise integration by parts delivers the following broken ultraweak variational formulation:

$$
(6.2)\quad \begin{cases} \text{Find } (\boldsymbol{u},\boldsymbol{\sigma},\widehat{\boldsymbol{u}},\widehat{\boldsymbol{t}}_n)\in\boldsymbol{L}^2(\Omega)\times\boldsymbol{L}^2(\Omega;\mathbb{S})\times\boldsymbol{H}_0^{1/2}(\mathcal{S})\times\boldsymbol{H}^{-1/2}(\mathcal{S}) : \\ \big(J^{-1}\mathsf{S}:\boldsymbol{J}\boldsymbol{\sigma},\overline{\boldsymbol{J}}\boldsymbol{\tau}\big)_\Omega + \big(\boldsymbol{u},\boldsymbol{\nabla}_h\cdot\boldsymbol{\tau}\big)_\Omega - \big\langle\widehat{\boldsymbol{u}},\boldsymbol{\tau}\cdot\boldsymbol{n}\big\rangle_{\partial\mathcal{T}} = 0\,, \qquad \forall\,\boldsymbol{\tau}\in\boldsymbol{H}(\mathrm{div},\mathcal{T};\mathbb{S})\,, \\ \qquad \big(\boldsymbol{\sigma},\boldsymbol{\nabla}_h\boldsymbol{v}\big)_\Omega - \big(J\rho\omega^2\boldsymbol{u},\boldsymbol{v}\big)_\Omega - \big\langle\widehat{\boldsymbol{t}}_n,\boldsymbol{v}\big\rangle_{\partial\mathcal{T}} = \big(\boldsymbol{f},\boldsymbol{v}\big)_\Omega\,, \qquad \forall\,\boldsymbol{v}\in\boldsymbol{H}^1(\mathcal{T})\,. \end{cases}
$$

Here, we have closed the problem with homogeneous displacement boundary conditions and introduced the mesh traction solution variable $\mathrm{tr}^{\partial\mathcal{T}}_{\mathrm{div}}\sim\boldsymbol{\sigma}\boldsymbol{n}$ [17, 21].

*Remark* 6.1. Denote $\mathfrak{v} = (\boldsymbol{\tau},\boldsymbol{v})$ and define $\mathcal{V} = \boldsymbol{H}(\mathrm{div},\mathcal{T};\mathbb{S})\times\boldsymbol{H}^1(\mathcal{T})$. The corresponding adjoint graph norm is simply

$$
\|\mathfrak{v}\|^2_{\mathcal{L}^*} = \|\,\boldsymbol{\nabla}_h\boldsymbol{v}-\overline{J}^{-1}\mathsf{S}_{\boldsymbol{J}}:\boldsymbol{\tau}\|^2_\Omega + \|\,\boldsymbol{\nabla}_h\cdot\boldsymbol{\tau}-\overline{J}\rho\omega^2\boldsymbol{v}\|^2_\Omega + \|\boldsymbol{\tau}\|^2_\Omega + \|\boldsymbol{v}\|^2_\Omega\,, \quad \forall\,\mathfrak{v}\in\mathcal{V}\,,
$$

where $[\mathsf{S}_{\boldsymbol{J}}]_{ijkl} = \overline{\boldsymbol{J}}_{ni}\overline{\boldsymbol{J}}_{mk}\mathsf{S}_{njml}$, where summation is implied over the repeated indices.

## 7. Numerical experiments I: set-up

We chose to verify the PML variational formulations (4.3), (5.2), and (6.2), derived above, by approximating well-known, unbounded-domain, closed-form exact solutions. Our approach, which is common in the literature, began by considering the two- and three-dimensional Green's functions corresponding to the model equations above. In this section, we state the Green's functions we used, define and illustrate the computation domains, including the PML, declare the boundary conditions used in our experiments, and define our stretching function for the PML.

7.1. **Green's functions.** Generally, Green's functions depend upon two positions in space, $\mathbf{x}, \mathbf{x}_0 \in \mathbb{R}^d$. From now on, we will simply assume that $\mathbf{x}_0 = \mathbf{0}$ because it greatly simplifies the exposition.

7.1.1. *Acoustics.* Let $\delta(\cdot)$ be the Dirac delta and $r(\mathbf{x}) = |\mathbf{x}|$ be the Euclidean norm of the coordinate $\mathbf{x} \in \mathbb{R}^d$. Now, consider the following second-order distributional Helmholtz equation in $\mathbb{R}^d$:

$$-\Delta p - \omega^2 p = \delta(r)\,,$$

with the Sommerfeld radiation boundary condition

$$\lim_{r\to\infty} r^{(d-1)/2}\left(\frac{\partial p}{\partial r} - i\omega p\right) = 0\,.$$

In two dimensions (i.e. $d = 2$), the solution to this problem is given by the Green's function

$$p^{\text{exact}} = \frac{i}{4}\mathcal{H}_0^{(1)}(\omega r)\,,$$

where $\mathcal{H}_0^{(1)}(\cdot)$ is the zeroth order Hankel function of the first kind [82]. In three dimensions (i.e. $d = 3$), the solution to this problem is given by the Green's function

$$p^{\text{exact}} = \frac{i}{4r}e^{i\omega r}\,.$$

7.1.2. *Electromagnetics.* From now on, we choose to only consider isotropic media and, therefore, we fix the electric permittivity $\epsilon = \epsilon_0$ and the magnetic permeability $\mu = \mu_0$. We will also set the material conductivity $\sigma = 0$.

Let $\boldsymbol{e} \in \mathbb{R}^d$ be a fixed unit vector, $|\boldsymbol{e}| = 1$, and let $\boldsymbol{e}_r(\mathbf{x}) \in \mathbb{R}^d$ be the corresponding radial unit vector in the $\mathbf{x}$-direction. Now, consider the following second-order distributional Maxwell equation in $\mathbb{R}^d$:

$$\frac{1}{\mu_0}\nabla\times\nabla\times\mathbf{E} - \omega^2\epsilon_0\mathbf{E} = i\omega\delta(r)\boldsymbol{e}\,,$$

with the Silver-Muller radiation boundary condition

$$\lim_{r\to\infty} r\Big(\nabla\times\mathbf{E} - ik_0\boldsymbol{e}_r\times(\boldsymbol{e}_r\times\mathbf{E})\Big) = 0\,.$$

Here, $k_0 = \omega\sqrt{\epsilon_0\mu_0}$ is the wavenumber. From now on, we will fix $\boldsymbol{e} = \boldsymbol{e}_x$, the constant unit vector in the $x$-direction.

In two dimensions (i.e. $d = 2$), the solution to this problem is given by the Green's function

$$\mathbf{E}^{\text{exact}} = \frac{i\omega\mu_0}{k_0^2}\begin{pmatrix} k_0^2 + \frac{\partial^2}{\partial x^2} \\ \frac{\partial^2}{\partial x\partial y}\end{pmatrix} g(r)\,, \qquad \text{where} \qquad g(r) = \frac{i}{4}\mathcal{H}_0^{(1)}(k_0 r)\,.$$

Likewise, in three dimensions (i.e. $d = 3$), the solution to this problem is given by the Green's function [83]

$$\mathbf{E}^{\text{exact}} = i\omega\mathbf{G}(r)\boldsymbol{e}_x\,, \qquad \text{where} \qquad \mathbf{G}(r) = \frac{\mu_0}{4r}\big(\mathbf{I} + k_0^{-2}\nabla\nabla^{\mathsf{T}}\big)e^{ik_0 r}\,.$$

7.1.3. *Elastodynamics.* In Section 6, the first-order systems (6.1) were expressed for arbitrary compliance tensors $\mathsf{S} = \mathsf{C}^{-1} : \mathbb{S} \to \mathbb{S}$. From now on, we choose to only consider isotropic homogeneous materials with uniform mass density $\rho = \rho_0$. Therefore, the elasticity tensor $\mathsf{C}$ is defined component-wise $\mathsf{C}_{ijkl} = \lambda\delta_{ij}\delta_{kl} + \mu(\delta_{ik}\delta_{jl} + \delta_{il}\delta_{jk})$, where $\lambda$ and $\mu$ are the Lamé parameters. Likewise, $\mathsf{S}_{ijkl} = \frac{1}{4\mu}(\delta_{ik}\delta_{jl} + \delta_{il}\delta_{jk}) - \frac{\lambda}{2\mu(3\lambda+2\mu)}\delta_{ij}\delta_{kl}$.

Now, consider the following second-order distributional linear elasticity equation in $\mathbb{R}^d$:

$$-\boldsymbol{\nabla}\cdot\big(\mathsf{C} : \boldsymbol{\nabla}\,\boldsymbol{u}\big) - \rho_0\omega^2\boldsymbol{u} = \delta(r)\boldsymbol{e}\,. \tag{7.1}$$

Define the compressive part of the wave field $\boldsymbol{u}_{\mathrm{p}} = -\frac{1}{k_{\mathrm{p}}}\nabla(\nabla\cdot\boldsymbol{u})$ and the shear part of the wave field $\boldsymbol{u}_{\mathrm{s}} = -\frac{1}{k_{\mathrm{s}}}\nabla\times(\nabla\times\boldsymbol{u})$, where $k_{\mathrm{p}} = \sqrt{\frac{\rho_0}{\lambda+2\mu}}\omega$ is the compressive wavenumber and $k_{\mathrm{s}} = \sqrt{\frac{\rho_0}{\mu}}\omega$ is the shear wavenumber. Under the isotropic constitutive law above, it can be shown that, $\boldsymbol{u} = \boldsymbol{u}_{\mathrm{p}} + \boldsymbol{u}_{\mathrm{s}}$. To close the unbounded problem (7.1), we impose the Kupradze-Sommerfeld radiation boundary conditions. Namely,

$$\lim_{r\to\infty} r^{(d-1)/2}\Big(\frac{\partial\boldsymbol{u}_{\mathrm{p}}}{\partial r} - ik_{\mathrm{p}}\boldsymbol{u}_{\mathrm{p}}\Big) = 0 \qquad \text{and} \qquad \lim_{r\to\infty} r^{(d-1)/2}\Big(\frac{\partial\boldsymbol{u}_{\mathrm{s}}}{\partial r} - ik_{\mathrm{s}}\boldsymbol{u}_{\mathrm{s}}\Big) = 0\,.$$

As in Section 7.1.2, let us consider a delta load in the $x$-direction, $\boldsymbol{e} = \boldsymbol{e}_x$. In two dimensions (i.e. $d = 2$), the two Cartesian components of the Green's function $\boldsymbol{u}^{\mathrm{exact}}$ are given by [82]

$$u_x^{\mathrm{exact}} = \frac{i}{4\mu}\Big(\Psi + \chi\frac{x^2}{r^2}\Big) \qquad \text{and} \qquad u_y^{\mathrm{exact}} = \frac{1}{4\mu}\Big(\chi\frac{xy}{r^2}\Big)\,,$$

where

$$\Psi = \mathcal{H}_0^{(1)}(k_{\mathrm{s}}r) + \Big(\frac{k_{\mathrm{p}}}{k_{\mathrm{s}}}\Big)^2\frac{\mathcal{H}_1^{(1)}(k_{\mathrm{p}}r)}{k_{\mathrm{p}}r} - \frac{\mathcal{H}_1^{(1)}(k_{\mathrm{s}}r)}{k_{\mathrm{s}}r} \qquad \text{and} \qquad \chi = \mathcal{H}_2^{(1)}(k_{\mathrm{s}}r) - \Big(\frac{k_{\mathrm{p}}}{k_{\mathrm{s}}}\Big)^2\mathcal{H}_2^{(1)}(k_{\mathrm{p}}r)\,.$$

Similarly, in three dimensions (i.e. $d = 3$), the three components of the Green's function $\boldsymbol{u}^{\mathrm{exact}}$ are given by

$$u_x^{\mathrm{exact}} = \frac{i}{4\pi\mu r}\Big(\Psi + \chi\frac{x^2}{r^2}\Big)\,, \qquad u_y^{\mathrm{exact}} = \frac{i}{4\pi\mu r}\Big(\chi\frac{xy}{r^2}\Big)\,, \qquad \text{and} \qquad u_z^{\mathrm{exact}} = \frac{1}{4\pi\mu r}\Big(\chi\frac{xz}{r^2}\Big)\,,$$

where

$$\Psi = \Big(\frac{k_{\mathrm{p}}}{k_{\mathrm{s}}}\Big)^2\Big(-\frac{i}{k_{\mathrm{p}}r} + \frac{1}{k_{\mathrm{p}}^2r^2}\Big)e^{ik_{\mathrm{p}}r} + \Big(1 + \frac{i}{k_{\mathrm{s}}r} - \frac{1}{k_{\mathrm{s}}^2r^2}\Big)e^{ik_{\mathrm{s}}r}$$

and

$$\chi = \Big(\frac{k_{\mathrm{p}}}{k_{\mathrm{s}}}\Big)^2\Big(1 + \frac{3i}{k_{\mathrm{p}}r} - \frac{3}{k_{\mathrm{p}}^2r^2}\Big)e^{ik_{\mathrm{p}}r} - \Big(1 + \frac{3i}{k_{\mathrm{s}}r} - \frac{3}{k_{\mathrm{s}}^2r^2}\Big)e^{ik_{\mathrm{s}}r}\,.$$

7.2. **Computational domains.** The computational domains which we have used to approximate the solutions of the unbounded problems above are given in Figure 7.1. Observe that, in each Cartesian direction, both meshes have eight elements in the interior subdomain and four elements in the PML region.

7.3. **Boundary conditions.** Tables 1 and 2 depict the boundary conditions applied to the two- and three-dimensional computational domains in Figure 7.1, respectively.

Note that each of the Green's functions above has a nontrivial singularity at the origin $\mathbf{x} = \mathbf{0}$ which has been tacitly avoided in the computational domains above. Therefore, in order to reproduce the Green's functions in computations, the values of the Green's functions were introduced into the boundary conditions on the (exterior) faces of the domain facing the singularity. Namely, in 2D, solution-consistent non-homogeneous boundary conditions were applied to $\{(x,y)\in\partial\Omega \mid x = 1 \vee y = 1\}$, and, in 3D, to $\{(x,y,z)\in\partial\Omega \mid x = 1 \vee y = 1 \vee z = 1\}$. On all other faces of the computational domain, homogeneous boundary conditions were applied. This has turned each boundary value problem above into a classical time-harmonic scattering problem.

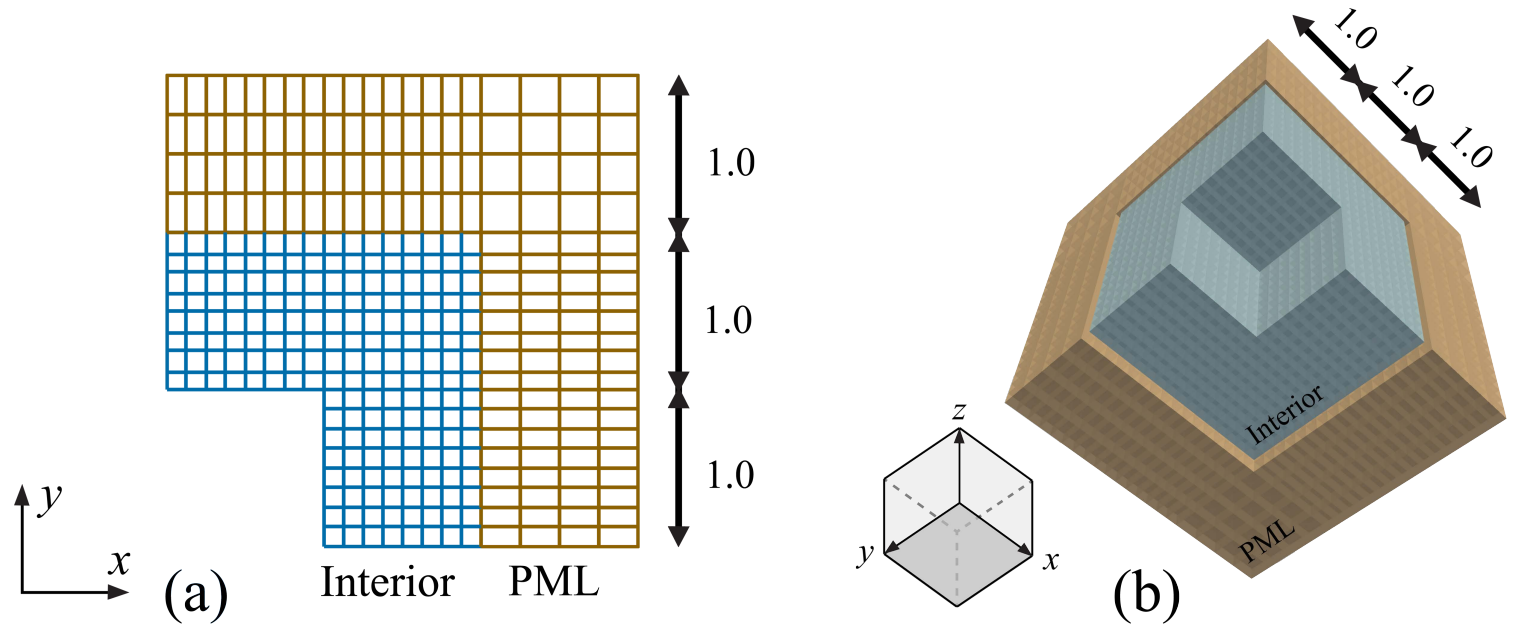

FIGURE 7.1. Computational domains and meshes separated into interior and PML subdomains: (a) 2D; (b) 3D. In the 3D illustration, a fictitious gap has been introduced between the interior and PML regions. This is for visualization purposes only, and is not present during computations.

| 2D problem | Boundary subsets. i.e. $\boldsymbol{x} = (x, y) \in \partial\Omega$ such that | | | |
|---|---|---|---|---|
| | $x = 1 \vee y = 1$ | $x = 0$ | $y = 0$ | $x = 3 \vee y = 3$ |
| Acoustics | $\widehat{p} = p^{\text{exact}}$ | $\widehat{\mathrm{u}}_n = 0$ | $\widehat{\mathrm{u}}_n = 0$ | $\widehat{p} = 0$ |
| Electromagnetics | $\widehat{\mathrm{E}} = \boldsymbol{n} \times \mathbf{E}^{\text{exact}}$ | $\widehat{\mathrm{E}} = 0$ | $\widehat{\mathrm{H}} = 0$ | $\widehat{\mathrm{E}} = 0$ |
| Elastodynamics | $\widehat{\boldsymbol{u}} = \boldsymbol{u}^{\text{exact}}$ | $\widehat{t}_{n_x} = \widehat{u}_y = 0$ | $\widehat{t}_{n_x} = \widehat{u}_y = 0$ | $\widehat{\boldsymbol{u}} = \mathbf{0}$ |

TABLE 1. Boundary conditions applied to the two-dimensional domain in Figure 7.1. In 2D, note that interface variables in (5.2), $\widehat{\mathrm{E}}$ and $\widehat{\mathrm{H}}$, are *scalar*-valued instead of vector-valued.

| 3D problem | Boundary subsets. i.e. $\boldsymbol{x} = (x, y, z) \in \partial\Omega$ such that | | | | |
|---|---|---|---|---|---|
| | $x = 1 \vee y = 1$ $\vee z = 1$ | $x = 0$ | $y = 0$ | $z = 0$ | $x = 3 \vee y = 3$ $\vee z = 3$ |
| Acoustics | $\widehat{p} = p^{\text{exact}}$ | $\widehat{\mathrm{u}}_n = 0$ | $\widehat{\mathrm{u}}_n = 0$ | $\widehat{\mathrm{u}}_n = 0$ | $\widehat{p} = 0$ |
| Electromagnetics | $\widehat{\mathbf{E}} = \boldsymbol{n} \times \mathbf{E}^{\text{exact}}$ | $\widehat{\mathbf{E}} = \mathbf{0}$ | $\widehat{\mathbf{H}} = 0$ | $\widehat{\mathbf{H}} = 0$ | $\widehat{\mathbf{E}} = \mathbf{0}$ |
| Elastodynamics | $\widehat{\boldsymbol{u}} = \boldsymbol{u}^{\text{exact}}$ | $\widehat{t}_{n_x} = \widehat{u}_y = \widehat{u}_z = 0$ | $\widehat{t}_{n_x} = \widehat{u}_y = \widehat{t}_{n_z} = 0$ | $\widehat{t}_{n_x} = \widehat{t}_{n_y} = \widehat{u}_z = 0$ | $\widehat{\boldsymbol{u}} = \mathbf{0}$ |

TABLE 2. Boundary conditions applied to the three-dimensional domain in Figure 7.1.

On the faces of the computational domain facing away from the singularity, we applied arbitrarily chosen homogeneous boundary conditions. Because each of these faces lie far within the PML, we expect the solution to decay rapidly enough that our choices here introduce negligible errors. On the final remaining faces — e.g. $\{(x, y) \in \partial\Omega \mid x = 0\}$, in 2D — the homogeneous boundary conditions were conveniently chosen to be compatible with the actual symmetries in the equations.

7.4. **Model parameters.** In all of our experiments, we used the fixed angular frequency $\omega = 6\pi$. In the computational domains we used (see Figure 7.1), this choice resulted in three acoustic wavelengths per domain. For electromagnetic model (5.1), we fixed $\epsilon_0 = \mu_0 = 1$ and for the elastodynamic model (6.1), we used the Lamé parameters $\lambda = 2$ and $\mu = \rho_0 = 1$. Therefore, there were also three electromagnetic wavelengths, three elastodynamic pressure wavelengths, and one and a half elastodynamic shear wavelengths, per domain.

7.5. **Stretching function.** Recall the stretching function given in (2.4). In both two and three dimensions, we used the explicit definition

$$\widetilde{x}_k = \begin{cases} x_k \,, & \text{if } 0 < x_k \leq l \,, \\ x_k + i\dfrac{C}{\omega}\left(\dfrac{x_k - l}{L - l}\right)^n, & \text{if } l < x_k \leq L \,, \end{cases}$$

where $l = 2$, $L = 3$, $n = 2$, and $C = 5$.

## 8. Numerical experiments II: implementation & results

8.1. **Discretization.** In this subsection, the symbol $p$ denotes an integer which corresponds to a given polynomial order.

Conforming finite element spaces for *non*-PML variational formulations corresponding to (4.3), (5.2), and (6.2) have been analyzed in several papers. For the acoustic wave equations, we choose to follow [4] and [7], for electromagnetics we follow [15], and for elastodynamics we follow [18]. Considering the important superconvergence discovery in [84], it is very likely that more efficient discretizations are possible.

Each of the non-superconvergent methods cited above employ a uniform polynomial order $p$ in the discretization of the trial space $\mathcal{U}$, which we respectively specify here to be $p = 4$. Meanwhile, each method also allows for a uniform polynomial $(p + \mathrm{d}p)$-order test space discretization to resolve the trial-to-test operator (3.2). In this article we always use $\mathrm{d}p = 1$.

8.2. **Software.** All of our computations were performed with the finite element software *hp*3D [85] which uses the ESEAS shape functions library to form the finite element spaces in each of the canonical energy spaces we require [86]. In each of our experiments non-homogeneous boundary conditions were applied to a non-trivial subset of the boundary. In order to apply these essential boundary conditions, we used projection-based interpolation [87]. In each of the energy spaces above, this is a fully-supported feature of the *hp*3D software.

8.3. **Results.** In our computational verification, for each model, in both two and three dimensions, less than 1% relative error was witnessed throughout the interior of the computational domains. The precise relative errors that our DPG methods delivered are recorded in Table 3. Explicitly, the relative error was

| Problem | 2D | 3D |
|---|---|---|
| Acoustics | 0.61286% | 0.52605% |
| Electromagnetics | 0.81578% | 0.73689% |
| Elastodynamics | 0.61861% | 0.58175% |

Table 3. The relative error in our numerical experiments, for each model problem, in both two and three dimensions.

Qualitatively, in each time-harmonic scattering model, in both two and three dimensions, throughout the entire computational domain, our experiments performed as expected. From a visual inspection of Figures 8.1–8.5, it is clear that each outgoing wave in the interior region steadily propagates from its source and then quickly dissipates within the PML region.

8.3.1. *Acoustics.* Figure 8.1 depicts the real part of the pressure variable $p$ from our DPG implementations of the two- and three-dimensional acoustic scattering problems described in Section 7.

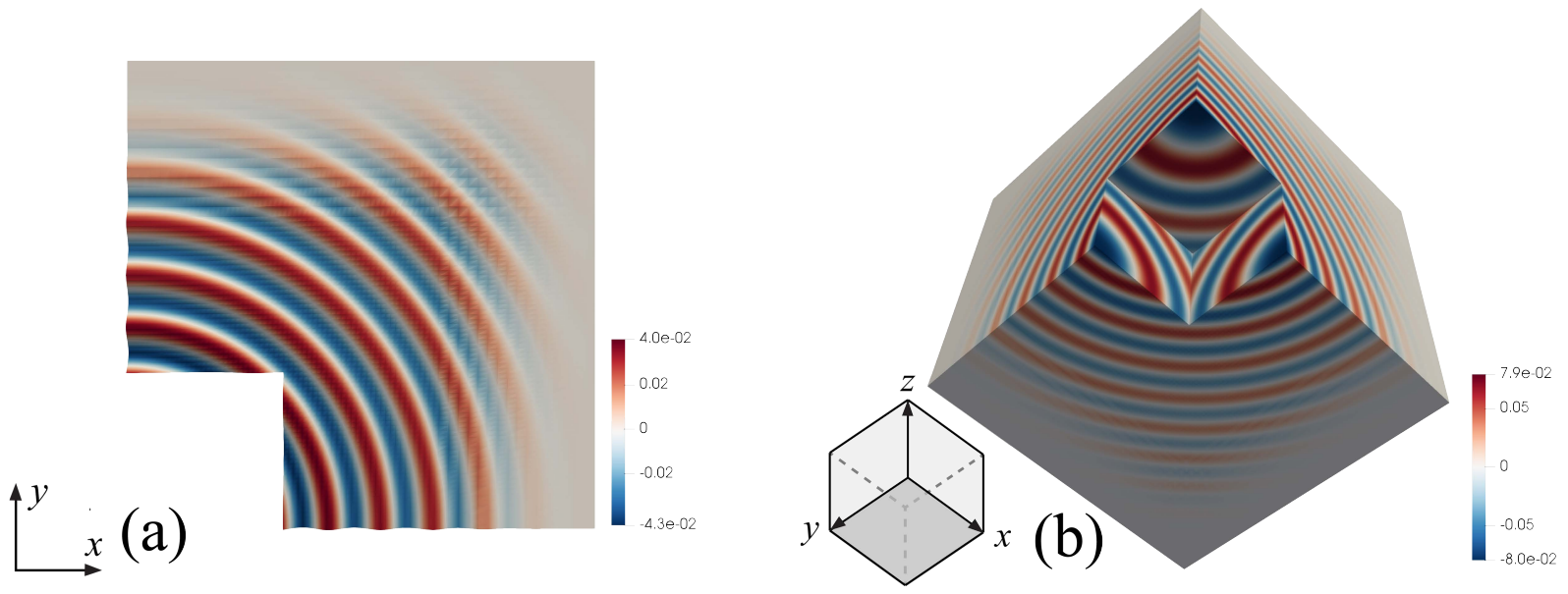

Figure 8.1. Time-harmonic acoustic wave scattering in the discrete pressure variable $p$ in: (a) 2D; (b) 3D. Here, only the real part of the solution is visualized.

8.3.2. *Electromagnetics.* Figures 8.2 and 8.3 depict the real part of the discrete electric field $\mathbf{E}$ from our DPG implementations of the two- and three-dimensional electromagnetic scattering problems described in Section 7.

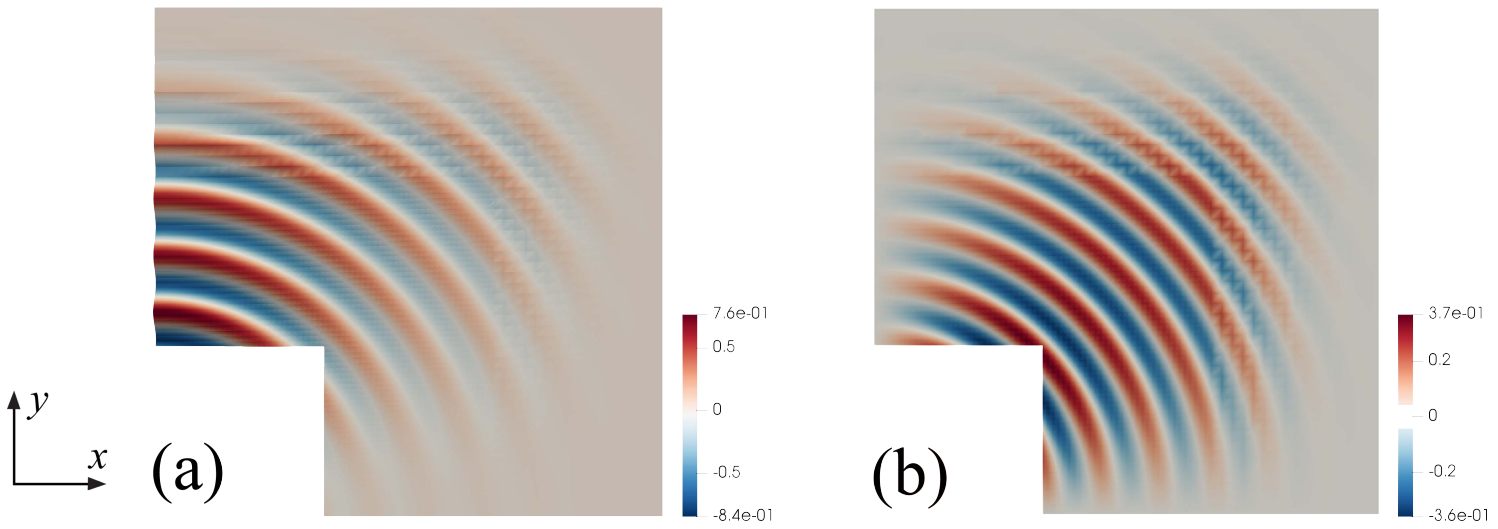

Figure 8.2. Time-harmonic electromagnetic wave scattering of the discrete electric field $\mathbf{E}$ in two dimensions: (a) the $x$-component of the electric field, $\mathrm{E}_x$; (b) the $y$-component of the electric field, $\mathrm{E}_y$. Here, only the real part of the solution is visualized.

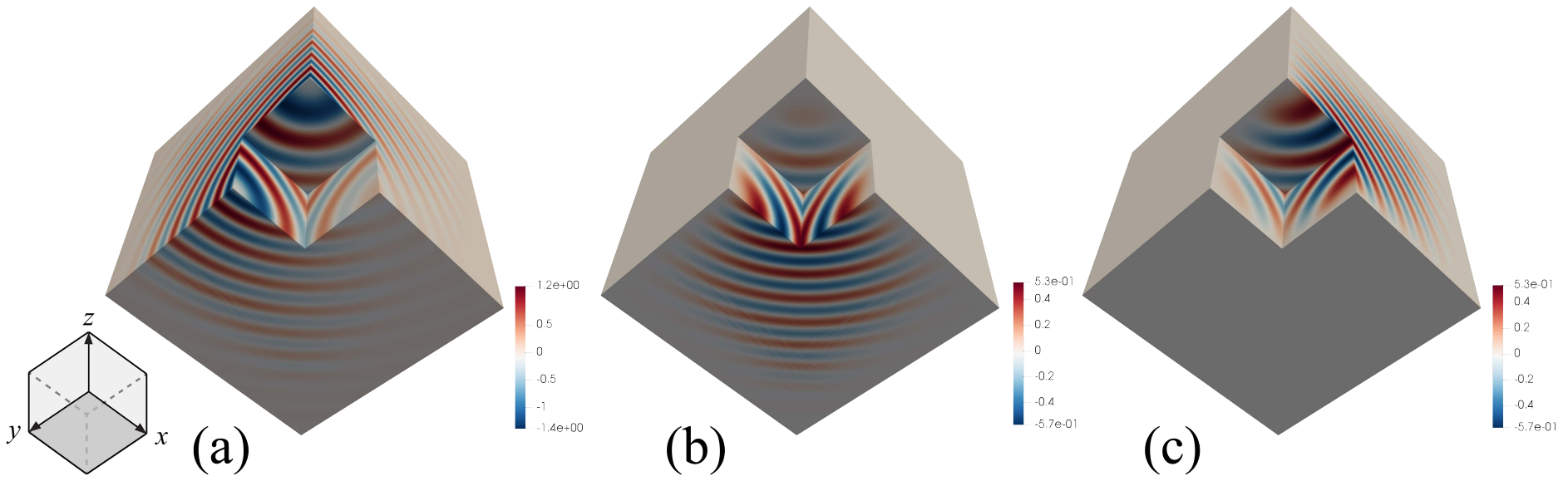

Figure 8.3. Time-harmonic electromagnetic wave scattering of the discrete electric field $\mathbf{E}$ in three dimensions: (a) the $x$-component of the electric field, $\mathrm{E}_x$; (b) the $y$-component of the electric field, $\mathrm{E}_y$; (c) the $z$-component of the electric field, $\mathrm{E}_z$. Again, only the real part of the solution is visualized.

8.3.3. *Elastodynamics.* Figures 8.4 and 8.5 depict the real part of the discrete displacement $\boldsymbol{u}$ from our DPG implementations of the two- and three-dimensional elastodynamic scattering problems described in Section 7.

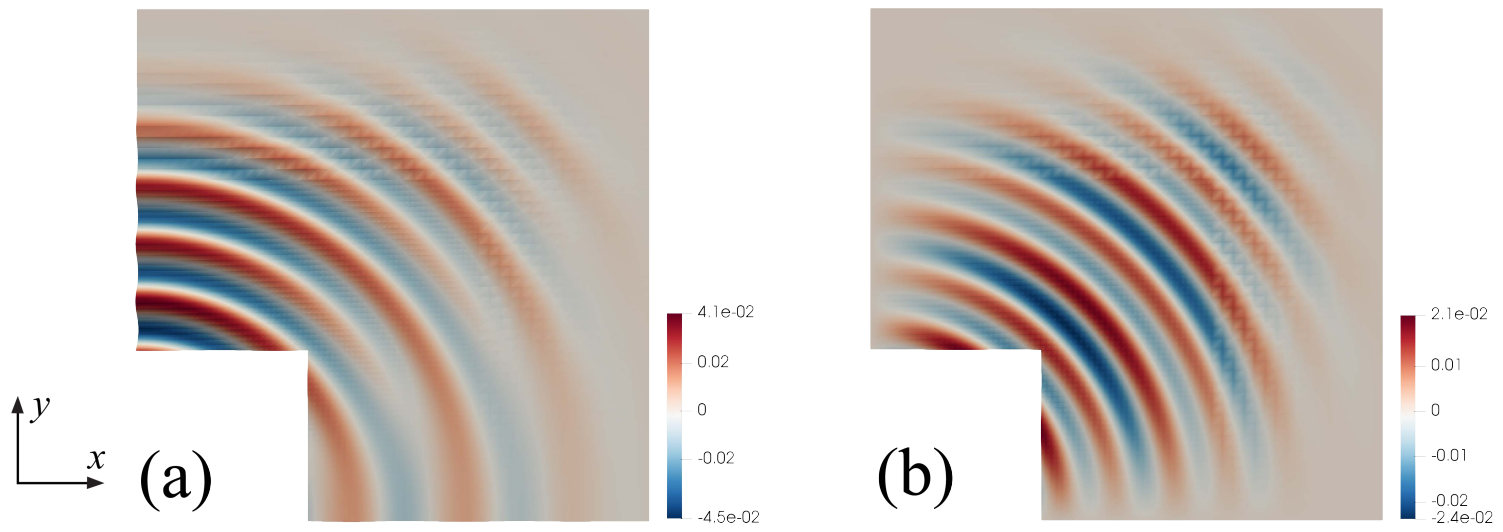

FIGURE 8.4. Time-harmonic elastodynamic wave scattering of the discrete displacement $\boldsymbol{u}$ in two dimensions: (a) the $x$-component of the displacement, $u_x$; (b) the $y$-component of the displacement, $u_y$. Here, only the real part of the solution is visualized.

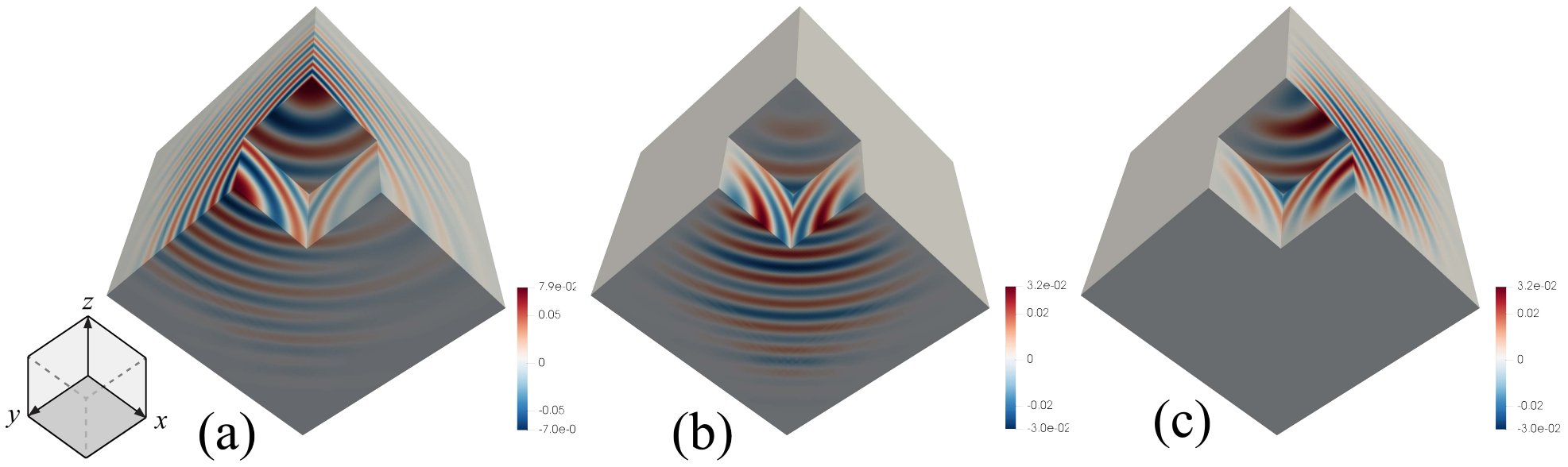

FIGURE 8.5. Time-harmonic elastodynamic wave scattering of the discrete displacement $\boldsymbol{u}$ in three dimensions: (a) the $x$-component of the displacement, $u_x$; (b) the $y$-component of the displacement, $u_y$; (b) the $z$-component of the displacement, $u_z$. Again, only the real part of the solution is visualized.

## 9. CONCLUSION

In this article, broken variational formulations are developed which treat time-harmonic problems posed on unbounded domains. Time-harmonic acoustics, electromagnetics, and elastodynamics models are considered in both two and three spatial dimensions. In the first case that is treated (acoustics), two different broken variational formulations with PML are eventually derived. The first derivation is new for DPG methods and the second was considered previously, in [18]. We remark that both are suitable for computation and that the two derivation strategies are known to deliver identical Bubnov-Galerkin methods. Eventually, the first strategy is judged to be most preferential and so it is then applied time-harmonic electromagnetics and elastodynamics models.

Numerical verification of the efficacy of the preferred PML variational formulations is reported on. The verification was performed in both two and three spatial dimensions, for each model, by comparing to the exact solutions — namely the Green's function — corresponding to each respective unbounded domains problem with the appropriate radiation boundary conditions. With the chosen PML and uniaxial stretching function, the numerical results are consistently positive across in all circumstances considered. Specifically, we report less than 1% relative error in each of our numerical verification experiments.

## References

[1] Moiola, A. and Spence, E. A., *Is the Helmholtz equation really sign-indefinite?*, SIAM Rev. 56.2 (2014), pp. 274–312.

[2] Demkowicz, L. and Gopalakrishnan, J., *Discontinuous Petrov–Galerkin (DPG) Method*, ICES Report 15-20, The University of Texas at Austin, 2015.

[3] Zitelli, J., Muga, I., Demkowicz, L., Gopalakrishnan, J., Pardo, D., and Calo, V. M., *A class of discontinuous Petrov–Galerkin methods. Part IV: The optimal test norm and time-harmonic wave propagation in 1D*, J. Comput. Phys. 230.7 (2011), pp. 2406–2432.

[4] Petrides, S. and Demkowicz, L. F., *An adaptive DPG method for high frequency time-harmonic wave propagation problems*, Comput. Math. Appl. 74.8 (2017), pp. 1999–2017.

[5] Ernst, O. G. and Gander, M. J., *Why it is difficult to solve Helmholtz problems with classical iterative methods?*, in: *Numerical Analysis of Multiscale Problems*, ed. by Graham, I. G., Hou, T. Y., Lakkis, O., and Scheichl, R., vol. 83, Springer Berlin Heidelberg, Berlin, Heidelberg, 2012, pp. 325–363.

[6] Li, X. and Xu, X., *Domain decomposition preconditioners for the discontinuous Petrov–Galerkin method*, ESAIM Math. Model. Num. 51.3 (2017), pp. 1021–1044.

[7] Keith, B., Petrides, S., Fuentes, F., and Demkowicz, L., *Discrete least-squares finite element methods*, Comput. Methods Appl. Mech. Engrg. 327 (2017), pp. 226–255.

[8] Demkowicz, L. and Gopalakrishnan, J., *A class of discontinuous Petrov–Galerkin methods. Part I: The transport equation*, Comput. Methods Appl. Mech. Engrg. 199.23-24 (2010), pp. 1558–1572.

[9] Demkowicz, L., Gopalakrishnan, J., Muga, I., and Zitelli, J., *Wavenumber explicit analysis of a DPG method for the multidimensional Helmholtz equation*, Comput. Methods Appl. Mech. Engrg. 213 (2012), pp. 126–138.

[10] Demkowicz, L. and Li, J., *Numerical simulations of cloaking problems using a DPG method*, Comput. Mech. 51.5 (2013), pp. 661–672.

[11] Bui-Thanh, T. and Ghattas, O., *A PDE-constrained optimization approach to the discontinuous Petrov–Galerkin method with a trust region inexact Newton-CG solver*, Comput. Methods Appl. Mech. Engrg. 278 (2014), pp. 20–40.

[12] Gopalakrishnan, J., Muga, I., and Olivares, N., *Dispersive and dissipative errors in the DPG method with scaled norms for Helmholtz equation*, SIAM J. Sci. Comput. 36.1 (2014), A20–A39.

[13] Gopalakrishnan, J. and Schöberl, J., *Degree and wavenumber [in]dependence of Schwarz preconditioner for the DPG method*, in: *Spectral and High Order Methods for Partial Differential Equations ICOSAHOM 2014*, ed. by Kirby, R., Berzins, M., and Hesthaven, J., vol. 106, Lecture Notes in Computational Science and Engineering, Springer, 2015, pp. 257–265.

[14] Gopalakrishnan, J. and Sepulveda, P., *A spacetime DPG method for acoustic waves*, ArXiv e-print arXiv:1709.08268 [math.NA] (2017).

[15] Carstensen, C., Demkowicz, L., and Gopalakrishnan, J., *Breaking spaces and forms for the DPG method and applications including Maxwell equations*, Comput. Math. Appl. 72.3 (2016), pp. 494–522,

[16] Niemi, A. H., Bramwell, J., and Demkowicz, L., *Discontinuous Petrov–Galerkin method with optimal test functions for thin-body problems in solid mechanics*, Comput. Methods Appl. Mech. Engrg. 200.9-12 (2011), pp. 1291–1300.

[17] Bramwell, J., Demkowicz, L., Gopalakrishnan, J., and Weifeng, Q., *A locking-free hp DPG method for linear elasticity with symmetric stresses*, Numer. Math. 122.4 (2012), pp. 671–707.

[18] Bramwell, J., *A discontinuous Petrov–Galerkin method for seismic tomography problems*, PhD thesis, The University of Texas at Austin, Austin, Texas, U.S.A., 2013.

[19] Gopalakrishnan, J. and Qiu, W., *An analysis of the practical DPG method*, Math. Comput. 83.286 (2014), pp. 537–552.

[20] Carstensen, C., Demkowicz, L., and Gopalakrishnan, J., *A posteriori error control for DPG methods*, SIAM J. Numer. Anal. 52.3 (2014), pp. 1335–1353.

[21] Keith, B., Fuentes, F., and Demkowicz, L., *The DPG methodology applied to different variational formulations of linear elasticity*, Comput. Methods Appl. Mech. Engrg. 309 (2016), pp. 579–609,

[22] Carstensen, C. and Hellwig, F., *Low-order discontinuous Petrov–Galerkin finite element methods for linear elasticity*, SIAM J. Numer. Anal. 54.6 (2016), pp. 3388–3410.

[23] Fuentes, F., Keith, B., Demkowicz, L., and Le Tallec, P., *Coupled variational formulations of linear elasticity and the DPG methodology*, J. Comput. Phys. 348 (2017), pp. 715–731.

[24] Fuentes, F., Demkowicz, L., and Wilder, A., *Using a DPG method to validate DMA experimental calibration of viscoelastic materials*, Comput. Methods Appl. Mech. Engrg. 325 (2017), pp. 748–765,

[25] Demkowicz, L., Gopalakrishnan, J., Nagaraj, S., and Sepulveda, P., *A spacetime DPG method for the Schrödinger equation*, SIAM J. Numer. Anal. 55.4 (2017), pp. 1740–1759.

[26] Bettess, P., *Infinite elements*, Int. J. Numer. Meth. Eng. 11.1 (1977), pp. 53–64.

[27] Demkowicz, L. and Shen, J., *A few new (?) facts about infinite elements*, Comput. Methods Appl. Mech. Engrg. 195.29 (2006), pp. 3572–3590.

[28] Keller, J. B. and Grote, M. J., *Exact nonreflecting boundary condition for elastic waves*, SIAM J. Numer. Anal. 60.3 (2000), pp. 803–819.

[29] Engquist, B. and Majda, A., *Absorbing boundary conditions for numerical simulation of waves*, Proc. Natl. Acad. Sci. U.S.A. 74.5 (1977), pp. 1765–1766.

[30] Engquist, B. and Majda, A., *Radiation boundary conditions for acoustic and elastic wave calculations*, Comm. Pure Appl. Math. 32.3 (1979), pp. 313–357.

[31] Bayliss, A. and Turkel, E., *Radiation boundary conditions for wave-like equations*, Comm. Pure Appl. Math. 33.6 (1980), pp. 707–725.

[32] Higdon, R. L., *Numerical absorbing boundary conditions for the wave equation*, Math. Comput. 49.179 (1987), pp. 65–90.

[33] Guddati, M. N. and Tassoulas, J. L., *Continued-fraction absorbing boundary conditions for the wave equation*, J. Comput. Acoust. 8.01 (2000), pp. 139–156.

[34] Berenger, J.-P., *A perfectly matched layer for the absorption of electromagnetic waves*, J. Comput. Phys. 114.2 (1994), pp. 185–200.

[35] Chew, W. C. and Weedon, W. H., *A 3D perfectly matched medium from modified Maxwell's equations with stretched coordinates*, Microw. Opt. Technol. Lett. 7.13 (1994), pp. 599–604.

[36] Chew, W. C., Jin, J. M., and Michielssen, E., *Complex coordinate stretching as a generalized absorbing boundary condition*, Microw. Opt. Technol. Lett. 15.6 (1997), pp. 363–369.

[37] Bramble, J. and Pasciak, J., *Analysis of a finite PML approximation for the three dimensional time-harmonic Maxwell and acoustic scattering problems*, Math. Comput. 76.258 (2007), pp. 597–614.

[38] Collino, F. and Monk, P., *The perfectly matched layer in curvilinear coordinates*, SIAM J. Sci. Comput. 19.6 (1998), pp. 2061–2090.

[39] Hohage, T., Schmidt, F., and Zschiedrich, L., *Solving time-harmonic scattering problems based on the pole condition II: convergence of the PML method*, SIAM J. Numer. Anal. 35.3 (2003), pp. 547–560.

[40] Lassas, M. and Somersalo, E., *On the existence and convergence of the solution of PML equations*, Computing 60.3 (1998), pp. 229–241.

[41] Turkel, E. and Yefet, A., *Absorbing PML boundary layers for wave-like equations*, Appl. Numer. Math. 27.4 (1998), pp. 533–557.

[42] Michler, C., Demkowicz, L., Kurtz, J., and Pardo, D., *Improving the performance of perfectly matched layers by means of hp-adaptivity*, Numer. Methods Partial Differ. Equ. 23.4 (2007), pp. 832–858.

[43] Vaziri Astaneh, A. and Guddati, M. N., *A two-level domain decomposition method with accurate interface conditions for the Helmholtz problem*, Int. J. Numer. Meth. Eng. 107.1 (2016), pp. 74–90.

[44] Bramble, J., Pasciak, J., and Trenev, D., *Analysis of a finite PML approximation to the three dimensional elastic wave scattering problem*, Math. Comput. 79.272 (2010), pp. 2079–2101.

[45] Chen, Z., Xiang, X., and Zhang, X., *Convergence of the PML method for elastic wave scattering problems*, Math. Comput. 85.302 (2016), pp. 2687–2714.

[46] Collino, F. and Tsogka, C., *Application of the perfectly matched absorbing layer model to the linear elastodynamic problem in anisotropic heterogeneous media*, Geophysics 66.1 (2001), pp. 294–307.

[47] Hastings, F. D., Schneider, J. B., and Broschat, S. L., *Application of the perfectly matched layer (PML) absorbing boundary condition to elastic wave propagation*, J. Acoust. Soc. Am. 100.5 (1996), pp. 3061–3069.

[48] Vaziri Astaneh, A., Urban, M. W., Aquino, W., Greenleaf, J. F., and Guddati, M. N., *Arterial waveguide model for shear wave elastography: implementation and in vitro validation*, Phys. Med. Biol. 62 (2017), pp. 5473–5494.

[49] Vaziri Astaneh, A. and Guddati, M. N., *Dispersion analysis of composite acousto-elastic waveguides*, Compos. Part B Eng. 130 (2017), pp. 200–216.

[50] Vaziri Astaneh, A. and Guddati, M. N., *Improved inversion algorithms for near-surface characterization*, Geophys. J. Int. 206.2 (2016), pp. 1410–1423.

[51] Bao, G. and Wu, H., *On the convergence of the solutions of PML equations for Maxwell's equations*, SIAM J. Numer. Anal 43 (2005), pp. 2121–2143.

[52] Teixeira, F. L. and Chew, W. C., *Advances in the theory of perfectly matched layers*, in: *Fast and Efficient Algorithms in Computational Electromagnetics*, ed. by Chew, W. C., Jin, J.-M., Michielssen, E., and Song, J., Artech House, Boston, 2001, chap. 7, pp. 283–346.

[53] Meza-Fajardo, K. C. and Papageorgiou, A. S., *A nonconvolutional, split-field, perfectly matched layer for wave propagation in isotropic and anisotropic elastic media: stability analysis*, Bull. Seismol. Soc. Am. 98.4 (2008), pp. 1811–1836.

[54] Bécache, E., Fauqueux, S., and Joly, P., *Stability of perfectly matched layers, group velocities and anisotropic waves*, J. Comput. Phys. 188.2 (2003), pp. 399–433.

[55] Appelö, D. and Kreiss, G., *A new absorbing layer for elastic waves*, J. Comput. Phys. 215.2 (2006), pp. 642–660.

[56] Duru, K. and Kreiss, G., *A well-posed and discretely stable perfectly matched layer for elastic wave equations in second order formulation*, Commun. Comput. Phys 11.5 (2012), pp. 1643–1672.

[57] Savadatti, S. and Guddati, M. N., *Absorbing boundary conditions for scalar waves in anisotropic media. Part 1: Time harmonic modeling*, J. Comput. Phys. 229.19 (2010), pp. 6696–6714.

[58] Savadatti, S. and Guddati, M. N., *Absorbing boundary conditions for scalar waves in anisotropic media. Part 2: Time-dependent modeling*, J. Comput. Phys. 229.18 (2010), pp. 6644–6662.

[59] Loh, P.-R., Oskooi, A. F., Ibanescu, M., Skorobogatiy, M., and Johnson, S. G., *Fundamental relation between phase and group velocity, and application to the failure of perfectly matched layers in backward-wave structures*, Phys. Rev. E 79 (2009), p. 065601,

[60] Druskin, V., Gütel, S., and Knizhnerman, L., *Near-optimal perfectly matched layers for indefinite helmholtz problems*, SIAM Rev. 58.1 (2016), pp. 90–116.

[61] Bermúdez, A., Hervella-Nieto, L., Prieto, A., and Rodríguez, R., *An optimal perfectly matched layer with unbounded absorbing function for time-harmonic acoustic scattering problems*, J. of Comput. Phys. 223.2 (2007), pp. 469–488.

[62] Collino, F. and Monk, P. B., *Optimizing the perfectly matched layer*, Comput. Methods Appl. Mech. Engrg. 164.1-2 (1998), pp. 157–171.

[63] Chew, W. and Jin, J., *Perfectly matched layers in the discretized space: an analysis and optimization*, Electromagnetics 16.4 (1996), pp. 325–340.

[64] Petropoulos, P. G., *An analytical study of the discrete perfectly matched layer for the time-domain Maxwell equations in cylindrical coordinates*, IEEE Trans. Antennas Propagat. 51.7 (2003), pp. 1671–1675.

[65] Pardo, D., Demkowicz, L., Torres-Verdín, C., and Michler, C., *PML enhanced with a self-adaptive goal-oriented hp-finite element method: Simulation of through-casing borehole resistivity measurements*, SIAM J. Sci. Comput. 30.6 (2008), pp. 2948–2964.

[66] Chen, Z. and Liu, X., *An adaptive perfectly matched layer technique for time-harmonic scattering problems*, SIAM J. Numer. Anal. 43.2 (2005), pp. 645–671.

[67] Chen, Z. and Wu, X., *An adaptive uniaxial perfectly matched layer method for time-harmonic scattering problems*, Numer. Math. Theor. Meth. Appl 1.2 (2008), pp. 113–137.

[68] Guddati, M., Lim, K., and Zahid, M., *Perfectly matched discrete layers for unbounded domain modeling*, in: *Computational Methods for Acoustics Problems*, ed. by F, M., Saxe-Coburg Publications, Scotland, 2008, pp. 69–98.

[69] Guddati, M. N. and Lim, K.-W., *Continued fraction absorbing boundary conditions for convex polygonal domains*, Int. J. Numer. Meth. Eng. 66.6 (2006), pp. 949–977.

[70] Vaziri Astaneh, A. and Guddati, M. N., *Efficient computation of dispersion curves for multilayered waveguides and half-spaces*, Comput. Methods Appl. Mech. Engrg. 300 (2016), pp. 27–46.

[71] Vaziriastaneh, A., *On the Forward and Inverse Computational Wave Propagation Problems.* PhD thesis, North Carolina State University, Raleigh, North Carolina, USA, 2016.

[72] Vaziri Astaneh, A. and Guddati, M. N., *WaveDisp: Dispersion Analysis Software for Immersed and Embedded Waveguides*, 2017, URL: http://www.WaveDisp.com.

[73] Vaziri Astaneh, A., Fuentes, F., Mora, J., and Demkowicz, L., *High-order polygonal discontinuous Petrov—Galerkin (PolyDPG) methods using ultraweak formulations*, Comput. Methods Appl. Mech. Engrg. 332 (2018), pp. 686–711.

[74] Matuszyk, P. J. and Demkowicz, L. F., *Parametric finite elements, exact sequences and perfectly matched layers*, Comput. Mech. 51.1 (2013), pp. 35–45.

[75] Demkowicz, L., *Various Variational Formulations and Closed Range Theorem*, ICES Report 15-03, The University of Texas at Austin, 2015.

[76] Gopalakrishnan, J., *Five lectures on DPG methods*, ArXiv e-prints arXiv:1306.0557 [math.NA] (2013).

[77] Roberts, N. V., *Camellia: a software framework for discontinuous Petrov–Galerkin methods*, Comput. Math. Appl. 68.11 (2014), pp. 1581–1604.

[78] Nagaraj, S., Petrides, S., and Demkowicz, L., *Construction of DPG Fortin operators for second order problems*, Comput. Math. Appl. 74.8 (2017), pp. 1964–1980.

[79] Bui-Thanh, T., Demkowicz, L., and Ghattas, O., *Constructively well-posed approximation methods with unity inf–sup and continuity constants for partial differential equations*, Math. Comput. 82.284 (2013), pp. 1923–1952.

[80] Keith, B., Demkowicz, L., and Gopalakrishnan, J., *DPG* method*, ArXiv e-prints arXiv:1710.05223 [math.NA] (2017).

[81] Keith, B., Vaziri Astaneh, A., and Demkowicz, L., *Goal-oriented adaptive mesh refinement for non-symmetric functional settings*, ArXiv e-prints arXiv:1711.01996 [math.NA] (2017).

[82] Kausel, E., *Fundamental solutions in elastodynamics: a compendium*, Cambridge University Press, 2006.

[83] Sarabandi, K., *Dyadic Green's function*, 2009, URL: `http://www.eecs.umich.edu/courses/eecs730/lect/DyadicGF_W09_port.pdf` (visited on 11/30/2017).

[84] Führer, T., *Superconvergence in the DPG method with ultra-weak formulation*, Comput. Math. Appl. 75.5 (2018), pp. 1705–1718.

[85] Demkowicz, L., Kurtz, J., Pardo, D., Paszyński, M., Rachowicz, W., and Zdunek, A., *Computing with hp Finite Elements. II. Frontiers: Three Dimensional Elliptic and Maxwell Problems with Applications*, Chapman & Hall/CRC, New York, 2007.

[86] Fuentes, F., Keith, B., Demkowicz, L., and Nagaraj, S., *Orientation embedded high order shape functions for the exact sequence elements of all shapes*, Comput. Math. Appl. 70.4 (2015), pp. 353–458,

[87] Demkowicz, L., *Polynomial exact sequences and projection-based interpolation with application to Maxwell equations*, in: *Mixed Finite Elements, Compatibility Conditions, and Applications*, ed. by Boffi, D. and Gastaldi, L., vol. 1939, Lecture Notes in Mathematics, Springer, Berlin, 2008, pp. 101–158.